\documentclass{amsart}

\usepackage{amsmath,amsthm,amsfonts,amscd,amssymb}

\newtheorem{thm}{Theorem}[section]
\newtheorem{cor}[thm]{Corollary}
\newtheorem{lem}[thm]{Lemma}
\newtheorem{prop}[thm]{Proposition}

\theoremstyle{definition}

\theoremstyle{remark}
\newtheorem{rem}{Remark}[section]

\begin{document}

\title{On distribution of well-rounded sublattices of $\mathbb Z^2$}
\author{Lenny Fukshansky}

\address{Department of Mathematics, Claremont McKenna College, 850 Columbia Avenue, Claremont, CA 91711-6420}
\email{lenny@cmc.edu}

\begin{abstract}
A lattice is called well-rounded if its minimal vectors span the corresponding Euclidean space. In this paper we completely describe well-rounded full-rank sublattices of ${\mathbb Z}^2$, as well as their determinant and minima sets. We show that the determinant set has positive density, deriving an explicit lower bound for it, while the minima set has density 0. We also produce formulas for the number of such lattices with a fixed determinant and with a fixed minimum. These formulas are related to the number of divisors of an integer in short intervals and to the number of its representations as a sum of two squares. We investigate the growth of the number of such lattices with a fixed determinant as the determinant grows, exhibiting some determinant sequences on which it is particularly large. To this end, we also study the behavior of the associated zeta function, comparing it to the Dedekind zeta function of Gaussian integers and to the Solomon zeta function of ${\mathbb Z}^2$. Our results extend automatically to well-rounded sublattices of any lattice $A {\mathbb Z}^2$, where $A$ is an element of the real orthogonal group $O_2({\mathbb R})$.
\end{abstract}

\maketitle
\tableofcontents

\def\A{{\mathcal A}}
\def\B{{\mathcal B}}
\def\C{{\mathcal C}}
\def\D{{\mathcal D}}
\def\E{{\mathcal E}}
\def\F{{\mathcal F}}
\def\x{{\mathcal H}}
\def\I{{\mathcal I}}
\def\J{{\mathcal J}}
\def\K{{\mathcal K}}
\def\L{{\mathcal L}}
\def\M{{\mathcal M}}
\def\Mm{{\mathfrak M}}
\def\Pp{{\mathfrak P}}
\def\Aa{{\mathfrak A}}
\def\N{{\mathcal N}}
\def\PP{{\mathcal P}}
\def\R{{\mathcal R}}
\def\s{{\mathcal S}}
\def\V{{\mathcal V}}
\def\W{{\mathcal W}}
\def\X{{\mathcal X}}
\def\Y{{\mathcal Y}}
\def\H{{\mathcal H}}
\def\OO{{\mathcal O}}
\def\aaa{{\mathbb A}}
\def\cee{{\mathbb C}}
\def\Nn{{\mathbb N}}
\def\pee{{\mathbb P}}
\def\que{{\mathbb Q}}
\def\real{{\mathbb R}}
\def\zed{{\mathbb Z}}
\def\gmn{{\mathbb G_m^N}}
\def\qbar{{\overline{\mathbb Q}}}
\def\DL{{\underline{\Delta}}}
\def\DU{{\overline{\Delta}}}
\def\eps{{\varepsilon}}
\def\vek{{\varepsilon_k}}
\def\ahat{{\hat \alpha}}
\def\bhat{{\hat \beta}}
\def\gt{{\tilde \gamma}}
\def\h{{\tfrac12}}
\def\ba{{\boldsymbol a}}
\def\be{{\boldsymbol e}}
\def\bei{{\boldsymbol e_i}}
\def\bc{{\boldsymbol c}}
\def\bm{{\boldsymbol m}}
\def\bk{{\boldsymbol k}}
\def\bi{{\boldsymbol i}}
\def\bl{{\boldsymbol l}}
\def\bq{{\boldsymbol q}}
\def\bu{{\boldsymbol u}}
\def\bt{{\boldsymbol t}}
\def\bs{{\boldsymbol s}}
\def\bv{{\boldsymbol v}}
\def\bw{{\boldsymbol w}}
\def\bx{{\boldsymbol x}}
\def\bX{{\boldsymbol X}}
\def\bz{{\boldsymbol z}}
\def\bwy{{\boldsymbol y}}
\def\bg{{\boldsymbol g}}
\def\bY{{\boldsymbol Y}}
\def\bL{{\boldsymbol L}}
\def\baa{{\boldsymbol\alpha}}
\def\bb{{\boldsymbol\beta}}
\def\bet{{\boldsymbol\eta}}
\def\bxi{{\boldsymbol\xi}}
\def\bo{{\boldsymbol 0}}
\def\bol{{\boldsymbol 1}_L}
\def\ep{\varepsilon}
\def\p{\boldsymbol\varphi}
\def\q{\boldsymbol\psi}
\def\WR{\operatorname{WR}}
\def\rank{\operatorname{rank}}
\def\aut{\operatorname{Aut}}
\def\lcm{\operatorname{lcm}}
\def\sgn{\operatorname{sgn}}
\def\spn{\operatorname{span}}
\def\md{\operatorname{mod}}
\def\Norm{\operatorname{Norm}}
\def\dim{\operatorname{dim}}
\def\det{\operatorname{det}}
\def\Vol{\operatorname{Vol}}
\def\rk{\operatorname{rk}}
\def\md{\operatorname{mod}}
\def\sqp{\operatorname{sqp}}

\section{Introduction and statement of results}

Let $N \geq 2$ be an integer, and let $\Lambda \subseteq \real^N$ be a lattice of full rank. Define the {\it minimum} of $\Lambda$ to be
$$|\Lambda| = \min_{\bx \in \Lambda \setminus \{\bo\}} \|\bx\|,$$
where $\|\ \|$ stands for the usual Euclidean norm on $\real^N$. Let
$$S(\Lambda) = \{ \bx \in \Lambda : \|\bx\| = |\Lambda| \}$$
be the set of {\it minimal vectors} of $\Lambda$. We say that $\Lambda$ is a {\it well-rounded} lattice (abbreviated WR) if $S(\Lambda)$ spans $\real^N$. WR lattices come up in a wide variety of different contexts, including sphere packing, covering, and kissing number problems, coding theory, and the linear Diophantine problem of Frobenius, just to name a few. Still, the WR condition is special enough so that one would expect WR lattices to be relatively sparce. However, in 2005 C. McMullen \cite{mcmullen} showed that in a certain sense {\it unimodular} WR lattices are ``well distributed'' among all {\it unimodular} lattices in $\real^N$, where a unimodular lattice is a lattice with determinant equal to 1.  More specifically, he proved the following theorem, from which he derived the 6-dimensional case of the famous Minkowski's conjecture for unimodular lattices.

\begin{thm} [\cite{mcmullen}] \label{mcmullen} Let $A \subseteq SL_N(\real)$ be the subgroup of diagonal matrices with positive diagonal entries, and let $\Lambda$ be a full-rank unimodular lattice in $\real^N$. If the closure of the orbit $A \Lambda$ is compact in the space of all full-rank unimodular lattices in $\real^N$, then it contains a WR lattice.
\end{thm}

\noindent
Notice that in a certain sense this is a statement about distribution of WR lattices in the space of all unimodular lattices in a fixed dimension. Motivated by this beautiful theorem, we want to investigate the distribution of WR sublattices of $\zed^N$, which is a natural arithmetic problem. For instance, for a fixed positive integer $t$, does there necessarily exist a WR subllatice $\Lambda \subseteq \zed^N$ so that $\det(\Lambda) = t$? If so, how many different such sublattices are there? The first trivial observation is that if $t = d^N$ for some $d \in \zed_{>0}$ and $I_N$ is the $N \times N$ identity matrix, then the lattice $\Lambda = (d I_N) \zed^N$ is WR with $\det(\Lambda) = t$ and $|\Lambda| = d$. It seems however quite difficult to describe {\it all} WR sublattices of $\zed^N$ in an arbitrary dimension $N$. This paper is concerned with providing such a description in dimension two.
\bigskip

From now on we will write $\WR(\Omega)$ for the set of all full-rank WR sublattices of a lattice $\Omega$; in this paper we will concentrate on $\WR(\zed^2)$. In section 3 we develop a certain parametrization of lattices in $\WR(\zed^2)$, which we then use to investigate the determinant set $\D$ of such lattices and to count the number of them for a fixed value of determinant. Specifically, let $\D$ be the set of all possible determinant values of lattices in $\WR(\zed^2)$, and let $\Mm$ be the set of all possible values of squared minima of these lattices, i.e. $\Mm = \{ |\Lambda|^2 : \Lambda \in \WR(\zed^2) \}$. It is easy to see that $\Mm$ is precisely the set of all positive integers, which are representable as a sum of two squares. Then it is interesting to understand how dense are these sets in $\zed_{>0}$. 

For any subset $\PP$ of $\zed$ and $M \in \zed_{>0}$, we write
$$\PP(M) = \{ n \in \PP : n \leq M\}.$$
Define {\it lower density} of $\PP$ in $\zed$ to be
$$\DL_{\PP} = \liminf_{M \rightarrow \infty} \frac{|\PP(M)|}{M},$$
and its {\it upper density} in $\zed$ to be
$$\DU_{\PP} = \limsup_{M \rightarrow \infty} \frac{|\PP(M)|}{M}.$$
Clearly, $0 \leq \DL_{\PP} \leq \DU_{\PP} \leq 1$. If $0 < \DL_{\PP}$, we say that $\PP$ {\it has density}, and if $\DL_{\PP} = \DU_{\PP}$, i.e. if $\lim_{M \rightarrow \infty} \frac{|\PP(M)|}{M}$ exists, we say that $\PP$ {\it has asymptotic density} equal to the value of this limit, which could be 0.

With this notation, we will show that $\D$ has density. More specifically, we prove the following.

\begin{thm} \label{dense} The determinant set $\D$ of lattices in $\WR(\zed^2)$ has representation
$$\D = \left\{ (a^2+b^2)cd\ :\ a,b \in \zed_{\geq 0},\ \max\{a,b\} >0,\ c,d \in \zed_{>0},\ 1 \leq \frac{c}{d} \leq \sqrt{3} \right\},$$
and lower density
\begin{equation}
\label{D_dens}
\DL_{\D} \geq \frac{3^{\frac{1}{4}}-1}{2 \cdot 3^{\frac{1}{4}}} \approx 0.12008216 \dots
\end{equation}
The minima set $\Mm$ has asymptotic density 0.
\end{thm}
\bigskip

We prove Theorem \ref{dense} in section 4. Now, if $\Lambda \in \WR(\zed^2)$, let $\bx,\bwy$ be a minimal basis for $\Lambda$, and let $\theta$ be the angle between the vectors $\bx$ and $\bwy$; it is a well known fact that in dimensions $\leq 4$ a lattice is always generated by vectors corresponding to its successive minima, so such a basis certainly exists (see, for instance, \cite{pohst}). Then there is a simple connection between the minimum and the determinant of $\Lambda$:
$$\det(\Lambda) = \|\bx\| \|\bwy\| \sin\theta = |\Lambda|^2 \sqrt{ 1 - \frac{\left( \bx^t \bwy \right)^2}{|\Lambda|^4}} = \sqrt{ |\Lambda|^4 - \left( \bx^t \bwy \right)^2 }.$$
Lemma \ref{gauss} below implies that $0 \leq |\bx^t \bwy| \leq \frac{|\Lambda|^2}{2}$. Therefore we have
$$\frac{\sqrt{3}\ |\Lambda|^2}{2} \leq \det(\Lambda) \leq |\Lambda|^2.$$
In view of this relation, it is especially interesting that the determinant set has positive density while the minima set has density 0. 
\bigskip
 
Next, for each $u \in \D$ we want to count the number of $\Lambda \in \WR(\zed^2)$ such that $\det(\Lambda) = u$. We need some additional notation. Suppose $t \in \zed_{>0}$ has prime factorization of the form
\begin{equation}
\label{prim_fact}
t = 2^w p_1^{2k_1} \dots p_s^{2k_s} q_1^{m_1} \dots q_r^{m_r},
\end{equation}
where $p_i \equiv 3\ (\md 4)$, $q_j \equiv 1\ (\md 4)$, $w \in \zed_{\geq 0}$, $k_i \in \frac{1}{2} \zed_{>0}$, and $m_j \in \zed_{>0}$ for all $1 \leq i \leq s$, $1 \leq j \leq r$. Let $\alpha(t)$ be the number of representations of $t$ as a sum of two squares ignoring order and signs, that is
\begin{equation}
\label{alpha}
\alpha(t) = \left| \left\{ (a,b) \in \zed^2_{\geq 0} : a^2+b^2 = t,\ a \leq b \right\} \right|.
\end{equation}
Also define
\begin{equation}
\label{alpha*}
\alpha_*(t) = \left| \left\{ (a,b) \in \zed^2_{\geq 0} : a^2+b^2 = t,\ a \leq b,\ \gcd(a,b)=1 \right\} \right|,
\end{equation}
for all $t > 2$, and define $\alpha_*(1) = \alpha_*(2) = \frac{1}{2}$. It is a well-known fact that $\alpha(t)$ is given by
\[ \alpha(t) = \left\{ \begin{array}{ll}
0 & \mbox{if any $k_i$ is a half-integer} \\
\frac{1}{2} B  & \mbox{if each $k_i$ is an integer and $B$ is even } \\
\frac{1}{2} \left( B - (-1)^w \right) & \mbox{if each $k_i$ is an integer and $B$ is odd,}
\end{array}
\right. \]
where $B=(m_1+1) \dots (m_r+1)$ (see, for instance \cite{silverman}, \cite{weisstein}). Clearly, when $t$ is squarefree, $\alpha_*(t)=\alpha(t)$. It is also a well-known fact that for $t$ as in (\ref{prim_fact})
\[ \alpha_*(t) = \left\{ \begin{array}{ll}
0 & \mbox{if $s \neq 0$ or $w > 1$} \\
2^{r-1} & \mbox{if $s=0$ and $w = 0$ or $1$.}
\end{array}
\right. \]
We also define the function
\begin{equation}
\label{beta}
\beta_{\nu}(t) = \left| \left\{ d \in \zed_{>0} : d\ |\ t\ \text{and } \frac{\sqrt{t}}{\nu} \leq d \leq \sqrt{t} \right\} \right|,
\end{equation}
for every $1 < \nu \leq 3^{1/4}$. The value of $\nu$ which will be particularly important to us is $\nu = 3^{1/4}$, therefore we define
$$\beta(t) = \beta_{3^{1/4}}(t).$$ 
We discuss the function $\beta_{\nu}(t)$ in more detail in section 2; at least it is clear that for each given $t$, $\beta_{\nu}(t)$ is effectively computable for every $\nu$. Finally, for any $t \in \zed_{>0}$ define
\[ \delta_1(t) = \left\{ \begin{array}{ll}
1 & \mbox{if $t$ is a square} \\
2 & \mbox{if $t$ is not a square,}
\end{array}
\right. \]
and
\[ \delta_2(t) = \left\{ \begin{array}{ll}
0 & \mbox{if $t$ is odd} \\
1 & \mbox{if $t$ is even, $\frac{t}{2}$ is a square} \\
2 & \mbox{if $t$ is even, $\frac{t}{2}$ is not a square.}
\end{array}
\right. \]
With this notation, we can state our second main theorem.

\begin{thm} \label{count} Let $u \in \zed_{>0}$, and let $\N(u)$ be the number of lattices in $\WR(\zed^2)$ with determinant equal to $u$. If $u=1$ or $2$, then $\N(u)=1$, the corresponding lattice being either $\zed^2$ or $\left( \begin{matrix} 1&-1 \\ 1&1 \end{matrix} \right) \zed^2$, respectively. Let $u>2$, and define
\[ t = t(u) = \left\{ \begin{array}{ll}
u & \mbox{if $u$ is odd} \\
\frac{u}{2}  & \mbox{if $u$ is even.}
\end{array}
\right. \]
Then:
\begin{eqnarray}
\label{N_formula}
\N(u) & = & \delta_1(t) \beta(t) + \delta_2(t) \beta \left(\frac{t}{2} \right) + 4 \mathop{\sum_{n|t, 1<n<t/2}}_{n\ \text{not a square}} \alpha_* \left(\frac{t}{n} \right) \beta(n) \nonumber \\
& + & 2 \mathop{\sum_{n|t, 1 \leq n<t/2}}_{n\ \text{a square}} \alpha_* \left(\frac{t}{n} \right) (2\beta(n)-1).
\end{eqnarray}
In particular, if $u \notin \D$, then the right hand side of (\ref{N_formula}) is equal to zero.
\end{thm}

Theorem \ref{count} can also be easily extended to a more general class of lattices. Namely, write $O_2(\real)$ for the real orthogonal group, then for every $A \in O_2(\real)$ and every $\bx,\bwy \in \real^2$ we have $(A\bx)^t (A\bwy) = \bx^t \bwy$, i.e. $O_2(\real)$ is the isometry group of $\real^2$ with respect to the Euclidean norm. Therefore, if $A \in O_2(\real)$ then $\Lambda \in \WR(A\zed^2)$ if and only if $A^t \Lambda \in \WR(\zed^2)$. This immediately implies the following result.

\begin{cor} \label{gen} Let $A \in O_2(\real)$. Then the determinant set and the minima set of lattices in $\WR(A\zed^2)$ are $\D$ and $\Mm$ respectively, as defined above. Moreover, for each $u \in \D$ the number of lattices in $\WR(A\zed^2)$ with determinant equal to $u$ is given by $\N(u)$ as in Theorem \ref{count}.
\end{cor}
\smallskip

We prove Theorem \ref{count} in section 5. In section 6 we use Theorem \ref{count} to work out simple examples of our formula in the case of prime power and product of two primes determinants. We also describe the ``orthogonal'' elements of $\WR(\zed^2)$, which come from ideals in Gaussian integers; these are quite sparse among all lattices in $\WR(\zed^2)$. We then derive easy to use bounds on $\N(u)$ and on the normal order of $\N(u)$. We also demonstrate examples of ``extremal'' sequences of determinant values, for which $\N(u)$ is especially large; see Corollary \ref{size_N}. In section 7 we derive a formula for the number of lattices in $\WR(\zed^2)$ of fixed minimum. 
\bigskip

In section 8 we study some basic properties of a zeta function, corresponding to the well-rounded sublattices of $\zed^2$. Namely, for $s \in \cee$ define
\begin{equation}
\label{WR_zeta}
\zeta_{\WR(\zed^2)}(s) = \sum_{\Lambda \in \WR(\zed^2)} (\det(\Lambda))^{-s} = \sum_{u=1}^{\infty} \N(u) u^{-s},
\end{equation}
where $\N(u)$ is as above. In particular, $\N(u) \neq 0$ if and only if $u \in \D$. For a Dirichlet series $\sum_{n=1}^{\infty} c_n n^{-s}$, we say that it has a {\it pole of order} $\mu$ at $s=s_0$, where $\mu$ and $s_0$ are positive real numbers, if
\begin{equation}
\label{pole_def}
0 < \lim_{s \rightarrow s_0^+} |s-s_0|^{\mu} \sum_{n=1}^{\infty} |c_n n^{-s}| < \infty.
\end{equation}
We will also say that such a Dirichlet series is {\it bounded from above (or below)} by a Dirichlet series $\sum_{n=1}^{\infty} b_n n^{-s}$, if $\sum_{n=1}^{\infty} |c_n n^{-s}| \leq \sum_{n=1}^{\infty} |b_n n^{-s}|$ (respectively, $\geq \sum_{n=1}^{\infty} |b_n n^{-s}|$). In section~8 we prove the following result.

\begin{thm} \label{zeta} Let the notation be as above, then $\zeta_{\WR(\zed^2)}(s)$ is analytic for all $s \in \cee$ with $\Re(s) > 1$, and is bounded from below by a Dirichlet series that has a pole of order 2 at $s=1$. Moreover, for every real $\eps >0$ there exists a Dirichlet series with a pole of order $2+\eps$ at $s=1$, which bounds $\zeta_{\WR(\zed^2)}(s)$ from above.
\end{thm}

\noindent
Notice that Theorem \ref{zeta} provides additional information about the growth of $\N(u)$. In section 8 we prove Theorem \ref{zeta} by means of considering the behavior of some related Dirichlet series, namely the generating functions of $\alpha_*$ and $\beta_{\nu}$. We should remark that we are not using the notion of a pole here in a sense that would imply the existence of an analytic continuation, but only to reflect on the growth of the coefficients; in fact, the arithmetic function $\N(u)$ behaves sufficiently erraticaly that one would doubt $\zeta_{\WR(\zed^2)}(s)$ having an analytic continuation to the left of $s=1$. We are now ready to proceed.
\bigskip

\section{A special divisor function}

As above, let $1 < \nu \leq 3^{1/4}$. In this section we briefly discuss bounds on the divisor function $\beta_{\nu}(t)$. Let
$$\M_{\nu}(t) = \left\{ d \in \zed_{>0} : d\ |\ t\ \text{and } \frac{\sqrt{t}}{\nu} \leq d \leq \sqrt{t} \right\}.$$

\begin{lem} \label{beta_gcd} If $d_1,d_2 \in \M_{\nu}(t)$, then $\gcd(d_1,d_2)>1$. 
\end{lem}

\proof
Suppose $d_1,d_2 \in \M_{\nu}(t)$, and $\gcd(d_1,d_2)=1$. Then $d_1d_2|t$, but
$$\frac{t}{\nu^2} < d_1d_2 \leq t.$$
Notice that $d_1d_2 \neq t$, since this would imply $d_1=d_2=\sqrt{t}$. Then
$$1 < \frac{t}{d_1d_2} < \nu^2 \leq \sqrt{3},$$
but $\frac{t}{d_1d_2} \in \zed$, which is a contradiction. 
\endproof

Lemma \ref{beta_gcd} implies in particular that $\M_{\nu}(t)$ can contain at most one prime $q$, and in this case every $d \in \M_{\nu}(t)$ must be divisible by $q$. Write $p(t)$ for the smallest prime divisor of $t$. Another immediate consequence of Lemma \ref{beta_gcd} is that gaps between two consecuitive elements of $\M_{\nu}(t)$ must be greater or equal than $p(t)$. Therefore, since $\beta_{\nu}(t)=|\M_{\nu}(t)|$, we obtain
\begin{equation}
\label{beta_bound1}
\beta_{\nu}(t) \leq \left[ \left( \frac{\nu - 1}{\nu p(t)} \right) \sqrt{t} \right] + 1,
\end{equation}
where $p(t) \geq 2$ for each $t \in \zed$, however for most $t$ better bounds are known.

Let us write $\tau(t)$ for the number of distinct divisors of $t$ and $\omega(t)$ for the number of distinct prime divisors of $t$ (see \cite{hall} for detailed information on $\tau(t)$ and $\omega(t)$). Hooley's $\Delta$-function of $t$ is defined by
$$\Delta(t) = \max_x \left| \left\{ d \in \zed_{>0} : d|t,\ e^x < d \leq e^{x+1} \right\} \right|.$$
If we take $x = \log \frac{\sqrt{t}}{\nu}$, then it is easy to see that
$$\M_{\nu}(t) \subseteq \left\{ d \in \zed_{>0} : d|t,\ e^x < d \leq e^{x+1} \right\},$$
and hence $\beta_{\nu}(t) \leq \Delta(t)$. Now, as stated in \cite{ten1} (see also \cite{hall}) a consequence of Sperner's theorem is that
\begin{equation}
\label{beta_bound2}
\beta_{\nu}(t) \leq \Delta(t) \leq O \left( \frac{\tau(t)}{\sqrt{\omega(t)}} \right),
\end{equation}
and if $t$ is squarefree the constant in $O$-notation is equal to 2. By Theorem 317 of \cite{hardy}, for any $\eps > 0$ there exists an integer $t_0(\eps)$ such that for all $t > t_0(\eps)$,
\begin{equation}
\label{tau_bound}
\tau(t) < 2^{(1+\eps) \frac{\log t}{\log \log t}}.
\end{equation}
Then combining (\ref{beta_bound2}) with (\ref{tau_bound}), we see that for any $\eps > 0$ there exists an integer $t_0(\eps)$ such that for all $t > t_0(\eps)$,\begin{equation}
\label{beta_bound3}
\beta_{\nu}(t) \leq O \left( t^{\frac{(1+\eps) \log 2}{\log \log t}} \right),
\end{equation}
which is much better than $O(\sqrt{t})$, the bound of (\ref{beta_bound1}), when $t$ is sufficiently large. In fact, for most $t$ we expect $\beta_{\nu}(t)$ to be even much smaller. A result of \cite{ten1} states that if $\psi(t)$ is any function such that $\psi(t) \rightarrow \infty$, as slowly as we wish, as $t \rightarrow \infty$, then
\begin{equation}
\label{beta_av}
\beta_{\nu}(t) \leq \Delta(t) < \psi(t) \log \log t
\end{equation}
for all positive integers $t$ in a sequence of asymptotic density 1. 
\bigskip

The function $\beta_{\nu}(t)$ also has a geometric interpretation. Write the prime decomposition for $t$ as
$$t = p_1^{e_1} \dots p_{\omega(t)}^{e_{\omega(t)}}$$
for corresponding distinct primes $p_1, \dots, p_{\omega(t)}$ and positive integers $e_1, \dots, e_{\omega(t)}$. If $d$ is a divisor of $t$, then
$$d = p_1^{x_1} \dots p_{\omega(t)}^{x_{\omega(t)}}$$
for some integers $0 \leq x_i \leq e_i$ for each $1 \leq i \leq \omega(t)$. Then $\frac{\sqrt{t}}{\nu} \leq d \leq \sqrt{t}$ if and only if
$$\log \frac{\sqrt{t}}{\nu} \leq \sum_{i=1}^{\omega(t)} x_i \log p_i \leq \log \sqrt{t}.$$
In other words, $\beta_{\nu}(t)$ is precisely the number of integer lattice points in the polytope $P_{\nu}(t)$ in $\real^{\omega(t)}$ bounded by the hyperplanes
\begin{equation}
\label{poly_pt}
x_i = 0,\ x_i = e_i\ \forall\ 1 \leq i \leq \omega(t),\  \sum_{i=1}^{\omega(t)} x_i = \log \frac{\sqrt{t}}{\nu},\ \sum_{i=1}^{\omega(t)} x_i = \log \sqrt{t}.
\end{equation}
In other words, $\beta_{\nu}(t) = \left|P_{\nu}(t) \cap \zed^{\omega(t)}\right|$, and $P_{\nu}(t)$ may be an irrational polytope. Counting integer lattice points in irrational polytopes is a very hard problem; a generating function for this problem is defined in \cite{barvinok}, but almost nothing seems to be known about it. 
\bigskip

\section{Parametrization of well-rounded lattices}

In this section we present an explicit description and a convenient parametrization of lattices in $\WR(\zed^2)$, which we later use to prove the main results of this paper. 

First we introduce some additional notation, following \cite{near:ort}. An ordered collection of linearly independent vectors $\{ \bx_1, \dots, \bx_k \} \subset \real^N$, $2 \leq k \leq N$, is called {\it nearly orthogonal} if for each $1 < i \leq k$ the angle between $\bx_i$ and the subspace of $\real^N$ spanned by $\bx_1, \dots, \bx_{i-1}$ is in the interval $\left[ \frac{\pi}{3}, \frac{2\pi}{3} \right]$. In other words, this condition means that for each $1 < i \leq k$
\begin{equation}
\label{near}
\frac{| \bx_i^t \bwy |}{\|\bx_i\| \|\bwy\|} \leq \frac{1}{2},
\end{equation}
for all non-zero vectors $\bwy \in \spn_{\real} \{ \bx_1, \dots, \bx_{i-1} \}$. The following result is Theorem 1 of \cite{near:ort}; in case $N=2$ this was proved by Gauss.

\begin{thm} [\cite{near:ort}] \label{no} Suppose that an ordered basis $\{ \bx_1, \dots, \bx_k \}$ for a lattice $\Lambda$ in $\real^N$ of rank $1 < k \leq N$ is nearly orthogonal. Then it contains a minimal vector of~$\Lambda$.
\end{thm}

\noindent
In particular, if all vectors $\bx_1, \dots, \bx_k$ of Theorem \ref{no} have the same norm, then $\Lambda$ is a WR lattice; we will call such a basis $\bx_1, \dots, \bx_k$ minimal.
\bigskip

For the rest of this paper we will restrict to the case $N=2$. Here is a first characterization of WR sublattices of $\zed^2$.

\begin{lem} \label{gauss} A sublattice $\Lambda \subseteq \zed^2$ of rank 2 is in $\WR(\zed^2)$ if and only if it has a basis $\bx,\bwy$ with
\begin{equation}
\label{gauss_cond}
\|\bx\| = \|\bwy\|,\ |\cos \theta| = \frac{| \bx^t \bwy |}{\|\bx\| \|\bwy\|} \leq \frac{1}{2},
\end{equation}
where $\theta$ is the angle between $\bx$ and $\bwy$. Moreover, if this is the case, then the set of minimal vectors $S(\Lambda) = \{ \pm \bx, \pm \bwy \}$. In particular, a minimal basis for $\Lambda$ is unique up to $\pm$ signs and reordering.
\end{lem}

\proof
Suppose first that $\Lambda$ contains a basis $\bx,\bwy$ satisfying (\ref{gauss_cond}). By Theorem \ref{no} this must be a minimal basis, meaning that $\Lambda$ is~WR.
\smallskip

Next assume that $\Lambda$ is WR, and let $\bx,\bwy \in S(\Lambda)$ be linearly independent vectors. It is a well known fact that for lattices of rank $\leq 4$ linearly independent minimal vectors form a basis, hence $\bx,\bwy$ is a basis for $\Lambda$, and $|\Lambda| = \|\bx\| = \|\bwy\|$. Let $\theta$ be the angle between $\bx$ and $\bwy$. We can assume without loss of generality that $\cos \theta > 0$: if not, replace $\bx$ with $-\bx$ or $\bwy$ with $-\bwy$. Notice that $\bo \neq \bx-\bwy \in \Lambda$, and
$$\|\bx - \bwy\| = \sqrt{ \|\bx\|^2 + \|\bwy\|^2 - 2 \bx^t \bwy } = |\Lambda| \sqrt{ 2 (1 - \cos \theta) }.$$
If $\cos \theta > \frac{1}{2}$, then $\|\bx - \bwy\| < |\Lambda|$, which is a contradiction. This proves (\ref{gauss_cond}), and also implies that the angle between two minimal linearly independent vectors in $\Lambda$ must lie in the interval $\left[ \frac{\pi}{3}, \frac{2\pi}{3} \right]$.
\smallskip

Now assume that $\Lambda \in \WR(\zed^2)$, and let $\bx,\bwy$ be a minimal basis for $\Lambda$, so $|\Lambda| = \|\bx\| = \|\bwy\|$. Let 
$$\theta_1,\theta_2,\theta_3,\theta_4$$
be angles between pairs of vectors $\{\bx,\bwy\}$, $\{\bwy,-\bx\}$, $\{-\bx,-\bwy\}$, and $\{-\bwy,\bx\}$ respectively. Since all of these vectors are in $\Lambda$ and have length $|\Lambda|$, it must be true that
$$\theta_1,\theta_2,\theta_3,\theta_4 \in \left[ \frac{\pi}{3}, \frac{2\pi}{3} \right].$$
On the other hand,
$$\theta_2 = \theta_4 = \pi - \theta_1,\ \theta_3 = \theta_1.$$
Assume there exists a vector $\bz \in \Lambda$ of length $|\Lambda|$ which is not equal to $\pm \bx, \pm \bwy$. Then all the angles it makes with the vectors $\pm \bx, \pm \bwy$ must lie in the interval $\left[ \frac{\pi}{3}, \frac{2\pi}{3} \right]$. This means that at least one of these angles must be equal to $\frac{\pi}{3}$, assume without loss of generality that this is the angle $\bz$ makes with $\bx$. Then
$$\bz = \left( \begin{matrix} z_1 \\ z_2 \end{matrix} \right) = \left( \begin{matrix} \cos \left( \frac{\pi}{3} \right)&-\sin \left( \frac{\pi}{3} \right) \\ \sin \left( \frac{\pi}{3} \right)&\cos \left( \frac{\pi}{3} \right) \end{matrix} \right) \left( \begin{matrix} x_1 \\ x_2 \end{matrix} \right) = \left( \begin{matrix} \frac{x_1}{2} - \frac{x_2\sqrt{3}}{2} \\ \frac{x_1\sqrt{3}}{2} + \frac{x_2}{2} \end{matrix} \right),$$
where $x_1,x_2 \in \zed$ are coordinates of $\bx$. Since $\bz \in \Lambda$, it must be true that $z_1,z_2 \in \zed$, but this is not possible. Hence a vector $\bz$ like this cannot exist, and this completes the proof.
\endproof

Next we develop a certain convenient explicit parametrization of lattices in $\WR(\zed^2)$. We start with lemmas describing two different families of such lattices.

\begin{lem} \label{abcd1} Let $a,b,c,d \in \zed$ be such that
\begin{equation}
\label{cond_esm1}
0 < |d| \leq |c| \leq \sqrt{3} |d|,\ \max \{|a|,|b|\} > 0.
\end{equation}
Then 
\begin{equation}
\label{esm1}
\Lambda = \left( \begin{matrix} ac+bd&ac-bd \\ bc-ad&bc+ad \end{matrix} \right) \zed^2
\end{equation}
is in $\WR(\zed^2)$ with 
\begin{equation}
\label{det_esm1}
\det(\Lambda) = 2(a^2+b^2)|cd|.
\end{equation}
\end{lem}

\proof
Suppose $a,b,c,d \in \zed$ satisfy (\ref{cond_esm1}). Let $\bx = \left( \begin{matrix} ac+bd \\ bc-ad \end{matrix} \right)$ and $\bwy = \left( \begin{matrix} ac-bd \\ bc+ad \end{matrix} \right)$, then
\begin{eqnarray}
\label{square}
\|\bx\|^2 & = & (ac+bd)^2 + (bc-ad)^2 \nonumber \\
& = & (a^2+b^2)(c^2+d^2) \nonumber \\
& = & (ac-bd)^2 + (bc+ad)^2 = \|\bwy\|^2.
\end{eqnarray}
Let $\Lambda = \spn_{\zed} \{\bx, \bwy\}$, then $\rk \Lambda = 2$. Let $\theta$ be the angle between $\bx$ and $\bwy$, and let $c = \gamma d$, where by (\ref{cond_esm1}), $1 \leq |\gamma| \leq \sqrt{3}$. Then, by (\ref{square}) and (\ref{cond_esm1})
\begin{eqnarray*}
|\cos(\theta)| & = & \frac{| \bx^t \bwy |}{\|\bx\| \|\bwy\|} = \frac{|(a^2+b^2)(c^2-d^2)|}{(a^2+b^2)(c^2+d^2)} \\
& = & \frac{c^2-d^2}{c^2+d^2} = \frac{\gamma^2-1}{\gamma^2+1} \leq \frac{1}{2}.
\end{eqnarray*}
Therefore $\theta \in \left[ \frac{\pi}{3}, \frac{2\pi}{3} \right]$, and so, by Lemma \ref{gauss}, $\Lambda$ is WR; (\ref{det_esm1}) follows. This completes the proof.
\endproof

\begin{lem} \label{abcd2} Let $a,b,c,d \in \zed$ be such that
\begin{equation}
\label{cond_esm2}
c^2+d^2 \geq 4|cd|,\ \max \{|a|,|b|\} > 0.
\end{equation}
Then 
\begin{equation}
\label{esm2}
\Lambda = \left( \begin{matrix} ac-bd&ad-bc \\ ad+bc&ac+bd \end{matrix} \right) \zed^2
\end{equation}
is in $\WR(\zed^2)$ with 
\begin{equation}
\label{det_esm2}
\det(\Lambda) = (a^2+b^2) |c^2-d^2|.
\end{equation}
\end{lem}

\proof
Suppose $a,b,c,d \in \zed$ satisfy (\ref{cond_esm2}). Let $\bx = \left( \begin{matrix} ac-bd \\ ad+bc \end{matrix} \right)$ and $\bwy = \left( \begin{matrix} ad-bc \\ ac+bd \end{matrix} \right)$, and define $\Lambda = \spn_{\zed} \{\bx, \bwy\}$. Then $\rk \Lambda = 2$, and $\|\bx\|, \|\bwy\|$ are the same as in (\ref{square}). Let $\theta$ be the angle between $\bx$ and $\bwy$. Then, by (\ref{square}) and (\ref{cond_esm2})
\begin{eqnarray*}
|\cos(\theta)| & = & \frac{| \bx^t \bwy |}{\|\bx\| \|\bwy\|} = \frac{2|cd|(a^2+b^2)}{(a^2+b^2)(c^2+d^2)} \\
& = & \frac{2|cd|}{c^2+d^2} \leq \frac{1}{2}.
\end{eqnarray*}
Therefore $\theta \in \left[ \frac{\pi}{3}, \frac{2\pi}{3} \right]$, and so, by Lemma \ref{gauss}, $\Lambda$ is WR; (\ref{det_esm2}) follows. This completes the proof.
\endproof

\begin{prop} \label{2-class} Suppose $\Lambda \in \WR(\zed^2)$. Then $\Lambda$ is either of the form as described in Lemma \ref{abcd1} or as in Lemma \ref{abcd2}.
\end{prop}

\proof
Suppose $\Lambda \in \WR(\zed^2)$, and let $\bx,\bwy \in \Lambda$ be a minimal basis
\begin{equation}
\label{c1}
\|\bx\|^2 = x_1^2 + x_2^2 = |\Lambda|^2 = y_1^2+y_2^2 = \|\bwy\|^2.
\end{equation}
Notice that due to (\ref{c1}) it must be true that either the pairs $x_1,y_1$ and $x_2,y_2$, or the pairs $x_1,y_2$ and $x_2,y_1$ are of the same parity. Indeed, suppose this is not true, then we can assume without loss of generality that $x_1,x_2$ are even and $y_1,y_2$ are odd. But then
$$x_1^2 + x_2^2 \equiv 0\ (\md\ 4),\ y_1^2+y_2^2 \equiv 2\ (\md\ 4),$$
which contradicts (\ref{c1}). Therefore, either
\begin{equation}
\label{c2}
\frac{x_1-y_1}{2},\ \frac{x_1+y_1}{2},\ \frac{y_2-x_2}{2},\ \frac{x_2+y_2}{2} \in \zed,
\end{equation}
or
\begin{equation}
\label{c3}
\frac{x_1-y_2}{2},\ \frac{x_1+y_2}{2},\ \frac{y_1-x_2}{2},\ \frac{x_2+y_1}{2} \in \zed.
\end{equation}
First assume (\ref{c2}) is true. Then let
\begin{equation}
\label{c4}
c = \gcd \left( \frac{x_1+y_1}{2}, \frac{x_2+y_2}{2} \right),\ a = \frac{x_1+y_1}{2c},\ b = \frac{x_2+y_2}{2c},\ d = \frac{(y_2-x_2)c}{x_1+y_1}.
\end{equation}
Clearly $a,b,c \in \zed$. We will now show that $d \in \zed$. Indeed,
$$d = \frac{(y_2-x_2)c}{x_1+y_1} = \frac{y^2_2-x^2_2}{(x_2+y_2)\left(\frac{x_1+y_1}{c}\right)},$$
and of course $(x_2+y_2)\ |\ (y^2_2-x^2_2)$. Also, by (\ref{c1})
$$\left(\frac{x_1+y_1}{c}\right)\ |\ (x_1+y_1)\ |\ (x_1^2-y_1^2) = (y^2_2-x^2_2),$$
and by definition of $c$ in (\ref{c4}),
$$\gcd\left( x_2+y_2,\ \frac{x_1+y_1}{c} \right) = 1,$$
which implies that
$$(x_2+y_2)\left(\frac{x_1+y_1}{c}\right)\ |\ (y^2_2-x^2_2),$$
and hence $d \in \zed$. With these definitions of $a,b,c,d$, it is easy to see that
$$x_1 = ac+bd,\ x_2 = bc-ad,\ y_1 = ac-bd,\ y_2 = bc+ad,$$
and hence $\Lambda$ is precisely of the form (\ref{esm1}). Moreover, since it is WR, Lemma \ref{gauss} implies that it must satisfy condition (\ref{gauss_cond}), which implies (\ref{cond_esm1}). This finishes the proof in case (\ref{c2}) is true. The proof in case (\ref{c3}) is true is completely analogous, in which case $\Lambda$ is of type (\ref{esm2}), and then (\ref{cond_esm2}) is satisfied.

\endproof
\smallskip

Suppose now that a lattice
$$\Lambda = \left( \begin{matrix} ac-bd&ad-bc \\ ad+bc&ac+bd \end{matrix} \right) \zed^2$$
with 
$$c^2+d^2 \geq 4|cd|,\ \max \{|a|,|b|\} > 0$$ 
as in Lemma \ref{abcd2}, and 
$$\det(\Lambda) = (a^2+b^2) |c^2-d^2|$$
is even. We will show that in this case $\Lambda$ can be represented in the form as in Lemma \ref{abcd1}. First assume that $a^2+b^2$ is even, then $a^2,b^2$, and hence $a,b$, must be of the same parity, meaning that $a+b$ and $a-b$ are even. Define
$$a_1 = \frac{a-b}{2},\ b_1 = \frac{a+b}{2},$$
then $a^2+b^2 = 2(a_1^2+b_1^2)$.
Let $c_1,d_1$ be such that $c_1 d_1 = c^2-d^2$, and
$$|c_1| = \max\{ |c-d|, |c+d| \},\ |d_1| = \min\{ |c-d|, |c+d| \}.$$
Suppose for instance that $c_1 = c+d$ and $d_1 = c-d$ (the argument is completely analogous in case $c_1 = c-d$ and $d_1 = c+d$). Then
$$c = \frac{c_1+d_1}{2},\ d = \frac{c_1-d_1}{2},$$
and so $4cd = c_1^2-d_1^2$. On the other hand $c^2+d^2 = \frac{c_1^2+d_1^2}{2}$. The fact that $c^2+d^2 \geq 4|cd|$ implies that
$$\frac{c_1^2+d_1^2}{2} \geq |c_1^2 - d_1^2| = c_1^2 - d_1^2,$$
since $|c_1| \geq |d_1|$, and so
$$|c_1| \leq \sqrt{3}\ |d_1|.$$
This choice of $a_1,b_1,c_1,d_1$ satisfies the conditions of (\ref{cond_esm1}), and it is easy to see that
\begin{equation}
\label{lat_eq}
\Lambda = \left( \begin{matrix} ac-bd&ad-bc \\ ad+bc&ac+bd \end{matrix} \right) \zed^2 = \left( \begin{matrix} a_1c_1+b_1d_1&a_1c_1-b_1d_1 \\ b_1c_1-a_1d_1&b_1c_1+a_1d_1 \end{matrix} \right) \zed^2,
\end{equation}
and $\det(\Lambda) = 2(a_1^2+b_1^2)|c_1d_1|$.
\smallskip

Next assume that $c^2-d^2$ is even. Then $(c+d)(c-d)$ is even and $(c+d)+(c-d) = 2c$ is even, which implies that $(c+d)$ and $(c-d)$ must both be even, in particular $c^2-d^2$ is divisible by $4$. Let $c_1,d_1$ be such that $4 c_1 d_1 = c^2-d^2$, and
$$|c_1| = \frac{1}{2} \max\{ |c-d|, |c+d| \},\ |d_1| = \frac{1}{2} \min\{ |c-d|, |c+d| \}.$$
By an argument as above, we can easily deduce again that
$$|d_1| \leq |c_1| \leq \sqrt{3}\ |d_1|.$$
Let 
$$a_1 = a-b,\ b_1 = a+b,$$
then $2(a^2+b^2) = a_1^2+b_1^2$, and so 
$$\det(\Lambda) = (a^2+b^2)|c^2-d^2| = 4(a^2+b^2)|c_1d_1| = 2(a_1^2+b_1^2)|c_1d_1|.$$
Once again, it is easy to check that with $a_1,b_1,c_1,d_1$ defined this way (\ref{lat_eq}) holds.
\bigskip

%%%%%%%%%%%%%%%%%%%%%%%%%%%%%
% In the argument above we were heavily relying on the basic fact that for every $a,b \in \zed$
% $$2(a^2+b^2) = a_1^2+b_1^2$$
% for some $a_1,b_1 \in \zed$. We also observed that for every $c,d \in \zed$ such that $c^2+d^2 \geq 4|cd|$, there exist $c_1,d_1 \in \zed$ with $|d_1| \leq |c_1| \leq \sqrt{3}\ |d_1|$ so that
% $$c^2-d^2 = c_1d_1.$$
%%%%%%%%%%%%%%%%%%%%%%%%%%%%

Let
\begin{equation}
\label{E1}
\E' = \left\{ (a,b,c,d) \in \zed^4\ :\ 0 < |d| \leq |c| \leq \sqrt{3} |d|,\ \max \{|a|,|b|\} > 0 \right\},
\end{equation}
and
\begin{eqnarray}
\label{O1}
\OO' = \{ (a,b,c,d) \in \zed^4 & : & 0 < |d| \leq |c| \leq \sqrt{3} |d|,\ \max \{|a|,|b|\} > 0, \nonumber \\
& & \text{so that } 2 \nmid (a^2+b^2) |cd| \}.
\end{eqnarray}
Define two classes of integral lattices 
\begin{equation}
\label{E}
\E = \left\{ \Lambda(a,b,c,d) = \left( \begin{matrix} ac+bd&ac-bd \\ bc-ad&bc+ad \end{matrix} \right) \zed^2 : (a,b,c,d) \in \E' \right\},
\end{equation}
and
\begin{equation}
\label{O}
\OO = \left\{ \Lambda(a,b,c,d) = \left( \begin{matrix} \frac{ac+ad+bd-bc}{2}&\frac{ac-ad-bc-bd}{2} \\ \frac{ac+bc+bd-ad}{2}&\frac{ac+ad+bc-bd}{2} \end{matrix} \right) \zed^2 : (a,b,c,d) \in \OO' \right\}.
\end{equation}
Then for every $\Lambda = \Lambda(a,b,c,d) \in \E$,
$$\det(\Lambda) = 2(a^2+b^2) |cd|,$$
is even, and for every $\Lambda = \Lambda(a,b,c,d) \in \OO$,
$$\det(\Lambda) = (a^2+b^2) |cd|$$
is odd. We proved the following theorem.

\begin{thm} \label{disjoint} The set $\WR(\zed^2)$ can be represented as the disjoint union of $\E$ and $\OO$. Moreover, the set of all possible determinants of lattices in $\WR(\zed^2) = \E \cup \OO$ is 
\begin{equation}
\label{dets}
\D = \{ (a^2+b^2) |cd| : a,b,c,d \in \zed,\ 0 < \max \{|a|,|b|\},\ 0 < |d| \leq |c| \leq \sqrt{3} |d| \}.
\end{equation}
\end{thm}

\begin{rem} \label{min_rem} Notice also that 
\[ |\Lambda|^2 = \left\{ \begin{array}{ll}
(a^2+b^2)(c^2+d^2) & \mbox{if $\Lambda \in \E$} \\
\frac{1}{2}(a^2+b^2)(c^2+d^2) & \mbox{if $\Lambda \in \OO$.}
\end{array}
\right. \]
Therefore the set of squared minima $\Mm(\E)$ of the lattices from $\E$ can be represented as
\begin{eqnarray}
\label{min_rep1}
\Mm(\E) = \{ (a^2+b^2)(c^2+d^2) & : & a,b,c,d \in \zed,\ 0 < \max \{|a|,|b|\},\nonumber \\
& & 0 < |d| \leq |c| \leq \sqrt{3} |d|,\ 2 | (a^2+b^2)cd \},
\end{eqnarray}
and the set of squared minima $\Mm(\OO)$ of the lattices from $\OO$ can be represented as
\begin{eqnarray}
\label{min_rep2}
\Mm(\OO) = \Big\{ \frac{1}{2} (a^2+b^2)(c^2+d^2) & : & a,b,c,d \in \zed,\ 0 < \max \{|a|,|b|\},\nonumber \\
& & 0 < |d| \leq |c| \leq \sqrt{3} |d|,\ 2 \nmid (a^2+b^2)cd \Big\}.
\end{eqnarray}
Then the set of squared minima $\Mm$ can be represented as $\Mm = \Mm(\E) \cup \Mm(\OO)$. 

%We will use this representation in section 6 to derive formulas for the number of such lattices with fixed minimum.
\end{rem}

\begin{cor} \label{monoid} The determinant set $\D$ in (\ref{dets}) and the squared minima set $\Mm$ are commutative monoids under multiplication.
\end{cor}

\proof
If $t_1 = (a_1^2+b_1^2)c_1d_1$ and $t_2 = (a_2^2+b_2^2)c_2d_2$ are in $\D$, then $t_1t_2 = (a_3^2+b_3^2)c_3d_3 \in \D$, where $a_3=a_1a_2+b_1b_2$, $b_3=b_1a_2-a_1b_2$, $c_3 = \pm \max\{|c_1d_2|,|d_1c_2|\}$, and $d_3 = \pm \min\{|c_1d_2|,|d_1c_2|\}$. It is also obvious that a product of two integers which are representable as sums of two squares is also representable as a sum of two squares. 
\endproof

Next we will use Theorem \ref{disjoint} to investigate the structure of the set $\D$ and to count the number of lattices in $\WR(\zed^2)$ of a fixed determinant.
\bigskip

\section{Proof of Theorem \ref{dense}}

The description of the set $\D$ in the statement of Theorem \ref{dense} follows immediately from (\ref{dets}). In this section we will mostly be concerned with deriving the estimate (\ref{D_dens}) for the lower density of $\D$. 

For each real number $1 < \nu \leq 3^{1/4}$, define the set
\begin{equation}
\label{prod_set}
\B_{\nu} = \left\{ n \in \zed_{>0} : \exists\ d \in \zed_{>0}\ \text{such that } d\ |\ n\ \text{and } \frac{\sqrt{n}}{\nu} \leq d \leq \sqrt{n} \right\}.
\end{equation}
Then notice that another description of the set $\D$ in (\ref{dets}) is
$$\D = \Mm \B_{3^{1/4}} = \{ m n : m \in \Mm,\ n \in \B_{3^{1/4}} \},$$
where $\Mm$ is the set of squared minima of lattices in $\WR(\zed^2)$, as before, so
\begin{equation}
\label{sum_set}
\Mm = \{ m \in \zed_{>0} : m = k^2+l^2\ \text{for some } k,l \in \zed \}.
\end{equation}
By a well known theorem of Fermat, the set $\Mm$ consists precisely of those positive integers $m$ in whose prime factorization every prime of the form $(4k+3)$ occurs an even number of times. Since $1 \in \Mm \cap \B_{3^{1/4}}$, we have $\Mm, \B_{3^{1/4}} \subset \D$; on the other hand, $6 \in \B_{3^{1/4}} \setminus \Mm$ and $2 \in \Mm \setminus \B_{3^{1/4}}$, hence $\Mm \subsetneq \D$ and $\B_{3^{1/4}} \subsetneq \D$. Moreover, $\D \subsetneq \zed_{>0}$, since for instance $3 \notin \D$. 

It is a well-known result of Landau (see, for instance \cite{motohashi}) that $\Mm$ has asymptotic density equal to 0, specifically
$$\lim_{M \rightarrow \infty} \frac{1}{M} \left| \{ n \in \Mm : n \leq M \} \right| = \lim_{M \rightarrow \infty} \frac{1}{\sqrt{\log M}} = 0.$$

Let us investigate the density of the sets $\B_{\nu}$ for a fixed $\nu \in (1,3^{1/4}]$. As before, for each $M \in \zed_{>0}$ we write
$$\B_{\nu}(M) = \{ n \in \B_{\nu} : n \leq M \}.$$
For each $n \in \zed_{>0}$ define
$$I_{\nu}(n) = \left\{n^2, n(n-1), \dots, n \left(n-\left[\left(\frac{\nu-1}{\nu}\right)n\right]\right)\right\}.$$
Notice that every $k \in \B_{\nu}(M)$ is of the form  $k = n(n-i)$ for some $n$ and $i \leq \left[\left(\frac{\nu-1}{\nu}\right)n\right]$, and so if $n \leq [\sqrt{M}]$, then $k \in \bigcup_{n=1}^{[\sqrt{M}]} I_{\nu}(n)$. For each $n \in \zed_{>0}$,
\begin{equation}
\label{I_card}
|I_{\nu}(n)| = \left[\left(\frac{\nu-1}{\nu}\right)n\right] + 1.
\end{equation}
There may also be some $k = n(n-i) \in \B_{\nu}(M)$ with $n > [\sqrt{M}]$ for some $i \leq \left(\frac{\nu-1}{\nu}\right)n$. Then $k = n^2 - ni \leq M$, and so $i \geq n - \frac{M}{n}$. It is easy to see that this is only possible if $n \leq \left[ \sqrt{\nu M} \right]$, and so for each $[\sqrt{M}] < n \leq \left[ \sqrt{\nu M} \right]$ define
$$J_{\nu,M}(n) = \left\{ n(n-i) : \left[ n - \frac{M}{n} \right] + 1 \leq i \leq \left[ \left(\frac{\nu-1}{\nu}\right)n \right] \right\}.$$
Clearly, $J_{\nu,M}(n) \subseteq I_{\nu}(n)$ for each such $n$, and
\begin{equation}
\label{B_contain_J}
\B_{\nu}(M) = \left( \bigcup_{n=1}^{[\sqrt{M}]} I_{\nu}(n) \right) \cup \left( \bigcup_{n=[\sqrt{M}]+1}^{\left[ \sqrt{\nu M} \right]} J_{\nu,M}(n)\right).
\end{equation}
For simplicity of approximation notice that
\begin{equation}
\label{B_contain}
\bigcup_{n=1}^{[\sqrt{M}]} I_{\nu}(n) \subseteq \B_{\nu}(M) \subseteq \bigcup_{n=1}^{\left[ \sqrt{\nu M} \right]} I_{\nu}(n).
\end{equation}
We immediately obtain an upper bound on $|\B_{\nu}(M)|$.

\begin{lem} \label{I_upp} For all $M \in \zed_{>0}$,
\begin{equation}
\label{I_upp1}
|\B_{\nu}(M)| \leq \frac{\nu-1}{2}\ M + \frac{\nu-1}{2 \sqrt{\nu}}\ \sqrt{M}.
\end{equation}
\end{lem}

\proof
Combining (\ref{B_contain}) and (\ref{I_card}), we obtain:
$$|\B_{\nu}(M)| \leq \sum_{n=1}^{\left[ \sqrt{\nu M} \right]} |I_{\nu}(n)| \leq  \frac{\nu-1}{\nu} \sum_{n=1}^{\left[ \sqrt{\nu M} \right]} n \leq \frac{\nu-1}{2 \nu} \sqrt{\nu M} \left( \sqrt{\nu M} + 1 \right).$$
The bound of (\ref{I_upp1}) follows.
\endproof

Next we want to produce a lower bound on $|\B_{\nu}(M)|$. For this we first consider the pairwise intersections of the sets $I_{\nu}(n)$.

\begin{lem} \label{I_int} Let $m < n \leq [\sqrt{M}]$.
\begin{trivlist}
\item (1)\ If $n > m \sqrt{\nu}$, or if $n \geq m \sqrt{\nu}$ and $\sqrt{\nu}$ is irrational, then $I_{\nu}(n) \cap I_{\nu}(m) = \emptyset$. 
\item (2)\ If $\gcd(m,n)=1$, then $I_{\nu}(n) \cap I_{\nu}(m) = \emptyset$.

\item (3)\ If $m < n < m \sqrt{\nu}$, then
$$|I_{\nu}(n) \cap I_{\nu}(m)| \leq \left[ \frac{\nu-1}{\nu}\ \gcd(m,n) \right] + 1.$$
\end{trivlist}
\end{lem}

\proof
Define $Q=\frac{\nu-1}{\nu}$. Let $m < n \leq [\sqrt{M}]$, and suppose that $k \in I_{\nu}(n) \cap I_{\nu}(m)$. Then
$$k = n(n-x) = m(m-y),$$
for some integers $0 \leq x \leq Qn$ and $0 \leq y \leq Qm$. Define a line 
$$L(n,m) = \{ (x,y) \in \real^2 : nx-my=n^2-m^2 \},$$
and a rectangular box
$$R(n,m) = \{ (x,y) \in \real^2 : 0 \leq x \leq Qn,\ 0 \leq y \leq Qm \}.$$
It follows immediately that
$$|I_{\nu}(n) \cap I_{\nu}(m)| = |L(n,m) \cap R(n,m) \cap \zed^2|.$$
The line $L(m,n)$ passes through the points $\left( \frac{n^2-m^2}{n}, 0 \right)$ and $\left( 0, -\frac{n^2-m^2}{n} \right)$, so in particular $L(n,m) \cap R(n,m) = \emptyset$ if $\frac{n^2-m^2}{n} > Qn$, i.e. if $n > m \sqrt{\nu}$. Also if $\sqrt{\nu}$ is irrational, then $m \sqrt{\nu}$ is never an integer, and so $I_{\nu}(n) \cap I_{\nu}(m) = \emptyset$ if $n \geq m \sqrt{\nu}$, proving (1).
\smallskip

Now suppose $m < n < m \sqrt{\nu}$, and let $(x,y) \in L(n,m) \cap R(n,m) \cap \zed^2$. Then 
$$y = -\frac{n(n-x)}{m} + m \in \zed_{>0},$$
hence $m\ |\ n(n-x)$. Clearly, $m \nmid n$, and so we must have
$$\lcm(m,n) = \frac{mn}{\gcd(m,n)}\ |\ n(n-x),$$
meaning that 
\begin{equation}
\label{gcd1}
\frac{m}{\gcd(m,n)}\ |\ n-x.
\end{equation}
In particular, if $\gcd(m,n)=1$, we must have $m\ |\ n-x$, but $n-x \leq n < 2m$, meaning that in order for $n-x$ to be divisible by $m$, it must be equal to $m$. This would imply that $x=n-m < \frac{n^2-m^2}{n}$, hence $y < 0$, meaning that $(x,y) \notin R(n,m)$, which is a contradiction. This proves (2).

Now assume that $(x_1,y_1),(x_1+t,y_2) \in L(n,m) \cap R(n,m) \cap \zed^2$, where $t$ is as small as possible. By (\ref{gcd1}), $\frac{m}{\gcd(m,n)}\ |\ n-x_1$ and $\frac{m}{\gcd(m,n)}\ |\ n-x_1-t$, so $\frac{m}{\gcd(m,n)}\ |\ t$, and by minimality of $t$ we must have $t = \frac{m}{\gcd(m,n)}$. Therefore $x$-coordinates of points in $L(n,m) \cap R(n,m) \cap \zed^2$ must satisfy 
$$\frac{n^2-m^2}{n} \leq x < Qn,$$
and the distance between $x$-coordinates of any two such points must be at least $\frac{m}{\gcd(m,n)}$. Hence,
\begin{eqnarray*}
|I_{\nu}(n) \cap I_{\nu}(m)| & = & |L(n,m) \cap R(n,m) \cap \zed^2| \leq \left[ \frac{Qn - \frac{n^2-m^2}{n}}{\frac{m}{\gcd(m,n)}} \right] + 1 \\
& = & \left[ \frac{(\nu m^2 - n^2) \gcd(m,n)}{\nu mn} \right] + 1 \leq \left[ \frac{\nu-1}{\nu}\ \gcd(m,n) \right] + 1.
\end{eqnarray*}
This proves (3).
\endproof

\begin{lem} \label{I_low} For all $M \in \zed_{>0}$,
\begin{equation}
\label{I_low1}
|\B_{\nu}(M)| > \frac{\nu-1}{2 \nu}\ M \left( 1 - \frac{ \log \log \sqrt{M}}{\log \sqrt{M}} \right)^2.
\end{equation}
\end{lem}

\proof
Let $N=\pi(\sqrt{M})$, i.e. the number of primes up to $\sqrt{M}$. It is a well-known fact that for all $\sqrt{M} \geq 11$,
\begin{equation}
\label{N_bound}
N \geq \frac{\sqrt{M}}{\log \sqrt{M}}.
\end{equation}
Hence suppose that $M \geq 121$, and let $p_1, \dots, p_N$ be all the primes up to $\sqrt{M}$ in ascending order. By part (2) of Lemma \ref{I_int},
$$I_{\nu}(p_i) \cap I_{\nu}(p_j) = \emptyset,$$
for all $1 \leq i \neq j \leq N$. As above, we let $Q=\frac{\nu-1}{\nu}$. Therefore, using (\ref{I_card}) we obtain:
\begin{equation}
\label{low2}
|\B_{\nu}(M)| \geq \sum_{i=1}^N |I_{\nu}(p_i)| \geq Q \sum_{i=1}^N p_i.
\end{equation}
A result of R. Jakimczuk \cite{jakimczuk} implies that
\begin{equation}
\label{prime_bound}
\sum_{i=1}^N p_i > \frac{N^2}{2}\ \log^2 N.
\end{equation}
The bound (\ref{I_low1}) follows upon combining (\ref{low2}) with (\ref{prime_bound}) and (\ref{N_bound}). 

Now assume that $M < 121$, so $\sqrt{M} < 11$. A direct verification shows that in this case $I_{\nu}(n) \cap I_{\nu}(m) = \emptyset$ for all $1 \leq n \neq m \leq [\sqrt{M}]$, and so
$$|\B_{\nu}(M)| \geq \sum_{n=1}^{[\sqrt{M}]} |I_{\nu}(n)| \geq Q \sum_{n=1}^{[\sqrt{M}]} n \geq \frac{\nu-1}{2 \nu}\ M.$$
This completes the proof.
\endproof

Combining Lemmas \ref{I_upp} and \ref{I_low}, we obtain the following result.

\begin{thm} \label{B_density}
$$\frac{\nu-1}{2 \nu} \leq \DL_{\B_{\nu}} \leq \DU_{\B_{\nu}} \leq \frac{\nu-1}{2}.$$
\end{thm}

\proof
Using (\ref{I_low1}), we see that
$$\DL_{\B_{\nu}} = \liminf_{M \rightarrow \infty} \frac{|\B_{\nu}(M)|}{M} \geq \frac{\nu-1}{2 \nu} \lim_{M \rightarrow \infty} \left( 1 - \frac{ \log \log \sqrt{M}}{\log \sqrt{M}} \right)^2 = \frac{\nu-1}{2 \nu},$$
and using (\ref{I_upp1}), we see that
$$\DU_{\B_{\nu}} = \limsup_{M \rightarrow \infty} \frac{|\B_{\nu}(M)|}{M} \leq \frac{\nu-1}{2} + \frac{\nu-1}{2 \sqrt{\nu}} \lim_{M \rightarrow \infty} \frac{1}{\sqrt{M}} = \frac{\nu-1}{2}.$$
This completes the proof.
\endproof

\begin{rem} It is also possible to produce bounds on $|\B_{\nu}(M)|$ using decomposition (\ref{B_contain_J}) instead of (\ref{B_contain}), and employing the fact that
$$|J_{\nu,M}(n)| = \left[ \left(\frac{\nu-1}{\nu}\right)n \right] - \left[ n - \frac{M}{n} \right].$$
It is also possible to employ the full power of Lemma \ref{I_int}, in particular part (3), to further refine the lower bound on $|\B_{\nu}(M)|$. These estimates however produce only marginally better constants, but much messier bounds in general.

We could also lift the restriction that $\nu \leq 3^{1/4}$ with essentially no changes to the arguments, but in any case the important situation is that with $\nu$ being close to 1, and we want to emphasize that the case of utmost importance to us is that with $\nu = 3^{1/4}$.
\end{rem}

Now (\ref{D_dens}) of Theorem \ref{dense} follows immediately by recalling that $\B_{3^{1/4}} \subseteq \D$, and applying Theorem \ref{B_density} with $\nu = 3^{1/4}$. The bounds on density of $\B_{\nu}$ are of independent interest, and will also be used in section 8 below to determine the order of the pole of the zeta function of well-rounded lattices.
\bigskip

\section{Proof of Theorem \ref{count}}

We start with a lemma which identifies all the 4-tuples from $\E' \cup \OO'$ which parametrize the same lattices.

\begin{lem} \label{transform1} Let $(a_1,b_1,c_1,d_1), (a_2,b_2,c_2,d_2) \in \E' \cup \OO'$. Then
\begin{equation}
\label{gamma1}
\Lambda(a_1,b_1,c_1,d_1) = \Lambda(a_2,b_2,c_2,d_2)
\end{equation}
if and only if there exists $0 \neq \gamma \in \que$ such that $(a_1,b_1,c_1,d_1)$ is equal to one of the following:
\begin{eqnarray}
\label{gamma2}
& & \left( \frac{a_2}{\gamma}, \frac{b_2}{\gamma}, \gamma c_2, \gamma d_2 \right),\ \left( \frac{b_2}{\gamma}, -\frac{a_2}{\gamma}, -\gamma d_2, \gamma c_2 \right),\ \left( -\frac{b_2}{\gamma}, \frac{a_2}{\gamma}, -\gamma d_2, \gamma c_2 \right), \nonumber \\
& & \left( -\frac{a_2}{\gamma}, -\frac{b_2}{\gamma}, \gamma c_2, \gamma d_2 \right),\ \left( -\frac{a_2}{\gamma}, -\frac{b_2}{\gamma}, -\gamma c_2, \gamma d_2 \right),\ \left( \frac{b_2}{\gamma}, -\frac{a_2}{\gamma}, \gamma d_2, \gamma c_2 \right), \nonumber \\
& & \left( -\frac{b_2}{\gamma}, \frac{a_2}{\gamma}, \gamma d_2, \gamma c_2 \right), \left( \frac{a_2}{\gamma}, \frac{b_2}{\gamma}, -\gamma c_2, \gamma d_2 \right).
\end{eqnarray}
\end{lem}

\proof
If $(a_1,b_1,c_1,d_1)$ is equal to one of the 4-tuples as in (\ref{gamma2}), then a direct verification shows that (\ref{gamma1}) is true. Suppose, on the other hand, that (\ref{gamma1}) is true. Let
$$\bx_1 = \left( \begin{matrix} a_1c_1+b_1d_1 \\ b_1c_1-a_1d_1 \end{matrix} \right),\ \bwy_1 = \left( \begin{matrix} a_1c_1-b_1d_1 \\ b_1c_1+a_1d_1 \end{matrix} \right),$$
and
$$\bx_2 = \left( \begin{matrix} a_2c_2+b_2d_2 \\ b_2c_2-a_2d_2 \end{matrix} \right),\ \bwy_2 = \left( \begin{matrix} a_2c_2-b_2d_2 \\ b_2c_2+a_2d_2 \end{matrix} \right),$$
if $(a_1,b_1,c_1,d_1), (a_2,b_2,c_2,d_2) \in \E'$, or
$$\bx_1 = \left( \begin{matrix} \frac{a_1c_1+a_1d_1+b_1d_1-b_1c_1}{2} \\ \frac{a_1c_1+b_1c_1+b_1d_1-a_1d_1}{2} \end{matrix} \right),\ \bwy_1 = \left( \begin{matrix} \frac{a_1c_1-a_1d_1-b_1c_1-b_1d_1}{2} \\ \frac{a_1c_1+a_1d_1+b_1c_1-b_1d_1}{2} \end{matrix} \right),$$
and
$$\bx_2 = \left( \begin{matrix} \frac{a_2c_2+a_2d_2+b_2d_2-b_2c_2}{2} \\ \frac{a_2c_2+b_2c_2+b_2d_2-a_2d_2}{2} \end{matrix} \right),\ \bwy_2 = \left( \begin{matrix} \frac{a_2c_2-a_2d_2-b_2c_2-b_2d_2}{2} \\ \frac{a_2c_2+a_2d_2+b_2c_2-b_2d_2}{2} \end{matrix} \right),$$
if $(a_1,b_1,c_1,d_1), (a_2,b_2,c_2,d_2) \in \OO'$. By Lemma \ref{gauss}, this means that the basis matrix $(\bx_1\ \bwy_1)$ for $\Lambda(a_1,b_1,c_1,d_1)$ must be equal to one of the following basis matrices for $\Lambda(a_2,b_2,c_2,d_2)$:
\begin{eqnarray*}
& & (\bx_2\ \bwy_2),\ (-\bx_2\ \bwy_2),\ (\bx_2\ -\bwy_2),\ (-\bx_2\ -\bwy_2), \\
& & (\bwy_2\ \bx_2),\ (-\bwy_2\ \bx_2),\ (\bwy_2\ -\bx_2),\ (-\bwy_2\ -\bx_2).
\end{eqnarray*}
A direct verification shows that in each of these cases $(a_1,b_1,c_1,d_1)$ is equal to one of the 4-tuples as in (\ref{gamma2}), in the same order. This completes the proof.
\endproof

\begin{rem} \label{sqr_free} Notice that in Lemma \ref{transform1} $(a_1,b_1,c_1,d_1)$ can be equal to the second, third, sixth, or seventh 4-tuple in (\ref{gamma2}) only if
$$|c_1|=|\gamma||d_2|,\ |d_1|=|\gamma||c_2|,$$
but on the other hand we know that $(a_1,b_1,c_1,d_1), (a_2,b_2,c_2,d_2) \in \E' \cup \OO'$. Combining these facts we obtain
\begin{equation}
\label{cd_eq}
|c_1| = |\gamma||d_2| \leq |\gamma||c_2| = |d_1| \leq |c_1|,
\end{equation}
which implies that there must be equality everywhere in (\ref{cd_eq}). In this case the determinant of the corresponding lattice $\Lambda$ is equal to $(a_1^2+b_1^2) c_1^2$ if $\Lambda \in \OO$ or to $2(a_1^2+b_1^2) c_1^2$ if $\Lambda \in \E$.
\end{rem}
\bigskip

\noindent
{\it Proof of Theorem \ref{count}.}
Let $u \in \D$. If $u=1,2$, the proof is by direct verification. Assume from here on that $u>2$, and let 
\[ t = t(u) = \left\{ \begin{array}{ll}
u & \mbox{if $u$ is odd} \\
\frac{u}{2}  & \mbox{if $u$ is even.}
\end{array}
\right. \]
Define
$$D(t) = \{ n \in \zed_{>0} : n | t \},$$
i.e. $D(t)$ is the set of positive divisors of $t$. Define
$$D_1(t) = \{ (c,d) \in D(t) \times D(t) : d \leq c \leq \sqrt{3} d,\ cd | t \}.$$
For each $(c,d) \in D_1(t)$, define
$$S_t(c,d) = \left\{ (a,b) \in \zed^2_{\geq 0} : a^2+b^2 = \frac{t}{cd},\ a \leq b \right\}.$$
Also let
$$T(t) = \{ (a,b,c,d) \in \zed^4_{\geq 0} : (c,d) \in D_1(t),\ (a,b) \in S_t(c,d) \}.$$
Define an equivalence relation on $T(t)$ by writing
$$(a_1,b_1,c_1,d_1) \sim (a_2,b_2,c_2,d_2)$$
if $(a_1,b_1,c_1,d_1) = \left( \frac{a_2}{\gamma},\frac{b_2}{\gamma},\gamma c_2,\gamma d_2 \right)$ for some $\gamma \in \que_{>0}$. Then let $T_1(t)$ be the set of all equivalence classes of elements of $T(t)$ under $\sim$, i.e. $T_1(t) = T(t)/\sim$. By abuse of notation, we will write $(a,b,c,d)$ for an element of $T_1(t)$. We first have the following lemma.

\begin{lem} \label{gcd_ab} For each equivalence class in $T_1(t)$ it is possible to select a unique representative $(a,b,c,d)$ with $\gcd(a,b)=1$. 
\end{lem}

\proof
Let $(a,b,c,d) \in T_1(t)$, and let $q=\gcd(a,b)$, then it is easy to see that
$$(a,b,c,d) \sim \left( \frac{a}{q},\frac{b}{q},qc,qd \right).$$
Moreover, suppose that $(a_1,b_1,c_1,d_1),(a_2,b_2,c_2,d_2) \in T_1(t)$ are such that 
$$\gcd(a_1,b_1)=\gcd(a_2,b_2)=1,$$ 
and
$$(a_1,b_1,c_1,d_1) \sim (a_2,b_2,c_2,d_2).$$
Then there exists $\gamma = \frac{s}{q} \in \que_{>0}$ with $\gcd(s,q)=1$ such that
$$a_1=\frac{q}{s} a_2,\ b_1=\frac{q}{s} b_2,\ c_1=\frac{s}{q} c_2,\ d_1=\frac{s}{q} d_2.$$
Then $s|a_2$, $s|b_2$, and so $s|\gcd(a_2,b_2)=1$, hence $s=1$. Also $q|qa_2=a_1$, $q|qb_2=b_1$, and so $q|\gcd(a_1,b_1)=1$, hence $q=1$. Therefore
$$(a_1,b_1,c_1,d_1) = (a_2,b_2,c_2,d_2).$$
This completes the proof.
\endproof

Therefore, for each $t$ we only need to count the lattices produced by the 4-tuples $(a,b,c,d) \in T_1(t)$ with $\gcd(a,b)=1$. Let $(a,b,c,d)$ be such a 4-tuple, then either $(a,b) = (0,1)$ or $a,b,c,d \neq 0$, since $\gcd(0,b)=b$. Moreover, $a \neq b$ unless $a=b=1$.

First assume $c \neq d$. If $(a,b) \neq (0,1),(1,1)$. By Lemma \ref{transform1} and Remark \ref{sqr_free} only the following 4-tuples produce the same lattice $\Lambda(a,b,c,d)$: 
\begin{eqnarray*}
& & (a,b,c,d),\ (-a,-b,-c,-d),\ (-a,-b,c,d),\ (a,b,-c,-d),\\
& & (a,b,-c,d),\ (a,b,c,-d),\ (-a,-b,-c,d),\ (-a,-b,c,-d).
\end{eqnarray*}
Then each $(a,b,c,d) \in T_1(t)$ gives rise to the four distinct lattices:
\begin{equation}
\label{sqrfr_lat1}
\Lambda(a,b,c,d),\ \Lambda(-a,b,c,d),\ \Lambda(b,a,c,d),\ \Lambda(-b,a,c,d), 
\end{equation}
since $\Lambda(a,b,d,c)=\Lambda(-b,a,c,d)$. Also, each of $(0,1,c,d),(1,1,c,d),(a,b,1,1) \in T_1(t)$ gives rise to the following pairs of distinct lattices, respectively:
\begin{eqnarray}
\label{sqrfr_lat2}
& & \Lambda(0,1,c,d),\ \Lambda(1,0,c,d); \nonumber \\
& & \Lambda(1,1,c,d),\ \Lambda(-1,1,c,d); \nonumber \\
& & \Lambda(a,b,1,1),\ \Lambda(-a,b,1,1).
\end{eqnarray}

Now suppose that $c=d$. Then 
$$\Lambda(a,b,c,c) = \Lambda(-b,a,c,c),\ \Lambda(-a,b,c,c) = \Lambda(b,a,c,c).$$
Hence, if $(a,b) \neq (0,1),(1,1)$, then each $(a,b,c,c) \in T_1(t)$ gives rise to two distinct lattices, $\Lambda(a,b,c,c)$ and $\Lambda(b,a,c,c)$. Notice also that
$$\Lambda(0,1,c,c)=\Lambda(1,0,c,c),\ \Lambda(1,1,c,c)=\Lambda(-1,1,c,c).$$
Hence 4-tuples $(0,1,c,c), (1,1,c,c) \in T_1(t)$ give rise to only one lattice each. 

\noindent
The formula for $\N(u)$, $u \in \D$, of Theorem \ref{count} follows. Also notice that if $u \in \zed_{>0} \setminus \D$, then for every divisor $n$ of $u$, either $\alpha_*(t/n)$ or $\beta(n)$ is equal to zero, and so the right hand side of (\ref{N_formula}) is equal to zero by construction.
This completes the proof of the theorem.
\boxed{ }
\smallskip

\begin{rem} \label{points} Our problem can be interpreted in terms of counting integral points on certain varieties. Let us say that two points 
$$\bx = (x_1,x_2,x_3,x_4)^t,\ \bwy  = (y_1,y_2,y_3,y_4)^t \in \real^4$$ 
are equivalent if there exists $U \in GL_2(\zed)$ such that
$$U \left( \begin{matrix} x_1&x_3 \\ x_2&x_4 \end{matrix} \right) = \left( \begin{matrix} y_1&y_3 \\ y_2&y_4 \end{matrix} \right).$$
Notice that the number of {\it all} full-rank sublattices of $\zed^2$ with determinant equal to $u$ is precisely the number of integral points on the hypersurface
$$x_1x_4 - x_2x_3 = u,$$
modulo this equivalence. This number is well known: one formula, for instance, is given by (\ref{all_sublattices}) below. On the other hand, by Lemma \ref{gauss}, the number of {\it well-rounded} full-rank sublattices of $\zed^2$ with determinant equal to $u$ is the number of integral points on the subset of the variety
$$x_1x_4 - x_2x_3 = u,\ x_1^2+x_2^2-x_3^2-x_4^2 = 0,$$
defined by the inequality
$$2|x_1x_3+x_2x_4| \leq x_1^2+x_2^2,$$
modulo the same equivalence. This makes direct counting much harder, and so our parametrization is quite useful.
\end{rem}
\bigskip

\section{Corollaries}

The first immediate consequence of Theorem \ref{count} is the following.

\begin{cor} \label{zero_even}  If $u \in \zed_{>0}$ is odd, then $\N(u) = \N(2u)$.
\end{cor}

\noindent
To demonstrate some examples of our formulas at work, we derive the following simpler looking expressions for the case of prime-power determinants.

\begin{cor} \label{primep} Let $p$ be a prime, $k \in \zed_{>0}$. Let $u = p^k$ or $2p^k$. Then
\[ \N(u) = \left\{ \begin{array}{ll}
0 & \mbox{if $p \equiv 3\ (\md 4)$ and $k$ is odd} \\
1 & \mbox{if $p \equiv 3\ (\md 4)$ and $k$ is even} \\
1 & \mbox{if $p=2$} \\
k+1 & \mbox{if $p \equiv 1\ (\md 4)$}
\end{array}
\right. \]
\end{cor}

\proof
First assume that $p \neq 2$. Define $t$ as in the statement of Theorem \ref{count}, then $t=p^k$. If $k$ is even, then by Theorem \ref{count}
\begin{eqnarray*}
\N(u) & = & \beta(p^k) + 4 \sum_{j=1}^{\frac{k}{2}} \alpha_*(p^{k+1-2j}) \beta(p^{2j-1}) + 2 \sum_{j=0}^{\frac{k}{2}-1} \alpha_*(p^{k-2j}) (2\beta(p^{2j})-1) \\
& = & 1 + 2 \sum_{j=0}^{\frac{k}{2}-1} \alpha_*(p^{k-2j}),
\end{eqnarray*}
since $\beta(p^{2j-1})=0$, and $\beta(p^{2j})=1$ for all $j$. If $p \equiv 3\ (\md 4)$, then $\alpha_*(p^{k-2j}) = 0$ for all $j$. If $p \equiv 1\ (\md 4)$, then $\alpha_*(p^{k-2j}) = 1$ for all $j$, in which case
$$\N(u) = 1 + 2 \sum_{j=0}^{\frac{k}{2}-1} 1 = k+1.$$

\noindent
Next assume $k>1$ is odd. Then, in the same manner as above,
\begin{eqnarray*}
\N(u) & = & 2 \beta(p^k) + 4 \sum_{j=1}^{\frac{k-1}{2}} \alpha_*(p^{k+1-2j}) \beta(p^{2j-1}) + 2 \sum_{j=0}^{\frac{k-1}{2}} \alpha_*(p^{k-2j}) (2\beta(p^{2j})-1) \\
& = & 2 \sum_{j=0}^{\frac{k-1}{2}} \alpha_*(p^{k-2j}),
\end{eqnarray*}
which is equal to 0 if $p \equiv 3\ (\md 4)$. If $p \equiv 1\ (\md 4)$, then
$$\N(u) = 2 \sum_{j=0}^{\frac{k-1}{2}} 1 = k+1.$$

\noindent
If $k=1$, then by Theorem \ref{count}
\[ \N(u) = 2 \alpha_*(p) = \left\{ \begin{array}{ll}
0 & \mbox{if $p \equiv 3\ (\md 4)$} \\
2 & \mbox{if $p \equiv 1\ (\md 4)$}.
\end{array}
\right. \]

Now assume that $p=2$, $u=p^k$, and $k>1$: the case $k=1$, i.e. $u=2$ is considered separately in the statement of Theorem \ref{count}. If $k$ is even, then
\begin{eqnarray*}
\N(u) & = & 2 \beta(2^{k-1}) + \beta(2^{k-2}) + 4 \sum_{j=1}^{\frac{k-2}{2}} \alpha_*(2^{k-2j}) \beta(2^{2j-1}) \\
& + & 2 \sum_{j=0}^{\frac{k-4}{2}} \alpha_*(2^{k-1-2j}) (2\beta(2^{2j})-1) = 1,
\end{eqnarray*}
since $\alpha_*(2^i)=0$ for all $i>1$, and $\beta(2^i)=0$ for all odd $i$. Now let $k$ be odd. Then
\begin{eqnarray*}
\N(u) & = & \beta(2^{k-1}) + 2 \beta(2^{k-2}) + 4 \sum_{j=1}^{\frac{k-3}{2}} \alpha_*(2^{k-2j}) \beta(2^{2j-1}) \\
& + & 2 \sum_{j=0}^{\frac{k-3}{2}} \alpha_*(2^{k-1-2j}) (2\beta(2^{2j})-1) = 1.
\end{eqnarray*}
This completes the proof. 
\endproof
\smallskip

\noindent
In precisely the same manner, we obtain the following formulas for the case when determinant is a product of two odd primes.

\begin{cor} \label{two_primes} If $u=p_1p_2$, where $p_1 < p_2$ are odd primes, then
\[ \N(u) = \left\{ \begin{array}{ll}
0 & \mbox{if $p_1$ or $p_2 \equiv 3\ (\md 4)$ and $p_2>\sqrt{3}p_1$} \\
2 & \mbox{if $p_1$ or $p_2 \equiv 3\ (\md 4)$ and $p_2 \leq \sqrt{3}p_1$} \\
4 & \mbox{if $p_1$ and $p_2 \equiv 1\ (\md 4)$ and $p_2>\sqrt{3}p_1$} \\
6 & \mbox{if $p_1$ and $p_2 \equiv 1\ (\md 4)$ and $p_2 \leq \sqrt{3}p_1$.}
\end{array}
\right. \]
\end{cor}

\proof
Direct verification.
\endproof

\noindent
The same way one can apply the formulas of Theorem \ref{count} to obtain explicit expressions for $\N(u)$ for many other instances of $u$ as well.
\smallskip

Notice that some of the lattices in $\WR(\zed^2)$ come from ideals in $\zed[i]$. Namely, let $u = a^2+b^2 \in \D$ and consider the lattices $\Lambda_1(a,b) = \left( \begin{matrix} a&-b \\ b&a \end{matrix} \right) \zed^2$ and $\Lambda_2(a,b) = \left( \begin{matrix} a&b \\ -b&a \end{matrix} \right) \zed^2$  with $\det(\Lambda_1) = \det(\Lambda_2) = u$. Let $I_1(a,b)$ and $I_2(a,b)$ be the ideals in $\zed[i]$ generated by $a+bi$ and $a-bi$ respectively, then $-b+ai = i(a+bi) \in I_1(a,b)$ and $b+ai = i(a-bi) \in I_2(a,b)$. Hence $I_1(a,b)$ and $I_2(a,b)$ map bijectively onto $\Lambda_1(a,b)$ and $\Lambda_2(a,b)$ respectively under the canonical mapping $x+iy \rightarrow \left( \begin{matrix} x \\ y \end{matrix} \right)$, and $\Lambda_1(a,b) = \Lambda_2(a,b)$ if and only if $b=0$, which can only happen when $u$ is a square. Notice that such representation is only possible for the determinant values $u$ which are also in the minima set $\Mm$; in other words, a full-rank WR sublattice of $\zed^2$ comes from an ideal in $\zed[i]$ if and only if it has an orthogonal basis. It is easy to see that the number of such lattices of determinant $u \in \Mm$, which is precisely the number of ideals of norm $u$ in $\zed[i]$, is equal to $2\alpha(u)$ if $u$ is not a square, and $2\alpha(u) + 1$ if $u$ is a square. With this in mind, we can now state the following immediate consequence of Corollaries \ref{primep} and \ref{two_primes}.

\begin{cor} \label{ideals} If $u \in \D$ is of the form $u = p^k, 2p^k$, where $p$ is a prime, or $u = p_1p_2$ where $\sqrt{3} p_1 < p_2$ are odd primes, then all lattices in $\WR(\zed^2)$ of determinant $u$ come from ideals of norm $u$ in $\zed[i]$.
\end{cor}

\noindent
This of course is not true in general, in fact the class of such lattices coming from ideals in $\zed[i]$ is quite thin. Notice in particular that in order for a lattice $\Lambda \in \WR(\zed^2)$ to come from an ideal of $\zed[i]$ it must first of all be true that $\det \Lambda \in \Mm$, which has density 0 versus the entire determinant set $\D$, which has positive density.
\smallskip

\begin{rem} Let $\Pp = \{p_1,p_2,\dots\}$ be the collection of all primes in the arithmetic progression $4n+1$. By Dirichlet's theorem on primes in arithmetic progressions, $\Pp$ is infinite. For each $p_i \in \Pp$ define $\PP_i = \{p_i^k,2p_i^k\}_{k=1}^{\infty}$. Then $\bigcup_{i=1}^{\infty} \PP_i \subset \D$, and Corollary \ref{primep} implies that for each $i$, 
$$\N(u) = \frac{\log u}{\log p_i} + 1,$$
for each $u \in \PP_i$. In other words, there are infinite sequences in $\D$ on which $\N(u)$ grows at least logarithmically in $u$. For comparison, it is a well known fact (see for instance \cite{bgruber}) that for any positive integer $u$ with prime factorization $u=q_1^{c_1} \dots q_m^{c_m}$ the number of {\it all} full-rank sublattices of $\zed^2$ with determinant $u$ is
\begin{equation}
\label{all_sublattices}
F(2,u) = \prod_{j=1}^m \frac{q_j^{c_j+1} - 1}{q_j - 1},
\end{equation}
which grows linearly in $u$. It is therefore interesting to exhibit sequences of determinant values $u$ for which $\N(u)$ is especially large.
\end{rem}
\smallskip

Recall that for an integer $u$, $\tau(u)$ and $\omega(u)$ are numbers of divisors and of prime divisors of $u$, respectively. We can report the following consequence of Theorem~\ref{count}. 

\begin{cor} \label{size_N} For each $u \in \zed_{>0}$,
\begin{equation}
\label{size_O1}
\N(u) \leq O \left( \tau(u)^2 2^{\omega(u)} \right) \leq O \left( \left( \frac{\sqrt{2} \log u}{\omega(u)} \right)^{2 \omega(u)} \right).
\end{equation}
Moreover,
\begin{equation}
\label{size_O2}
\N(u) < O \left( (\log u)^{\log 8} \right), 
\end{equation}
for all $u \in \D$ outside of a subset of asymptotic density 0. However, there exist infinite sequences $\{ u_k \}_{k=1}^{\infty} \subset \D$ such that for every $k \geq 1$
\begin{equation}
\label{size_O3}
\N(u_k) \geq (\log u_k)^k.
\end{equation}
For instance, there exists such a sequence with $u_k \leq \exp \left( O (k (\log k)^2) \right)$ and $\omega(u_k) = O(k \log k)$.
\end{cor}

\proof
Notice that the right hand side of (\ref{N_formula}) is the sum of at most $\tau(u)$ nonzero terms. Combining (\ref{beta_bound2}) with the formula for $\alpha_*$ in section 1, it follows that each of these terms is at most $O \left( \tau(u) 2^{\omega(u)} \right)$. Let $u$ have a prime decomposition of the form $u = p_1^{e_1} \dots p_n^{e_n}$, so $\omega(u)=n$, then:
$$\frac{\log u}{n} = \frac{1}{n} \sum_{i=1}^n e_i \log p_i \geq \left( \prod_{i=1}^n e_i \log p_i \right)^{\frac{1}{n}} \geq \left( O \left( \prod_{i=1}^n (e_i+1) \right) \right)^{\frac{1}{n}} = \left( O(\tau(u)) \right)^{\frac{1}{n}}.$$
This proves (\ref{size_O1}), and (\ref{size_O2}) follows from (\ref{size_O1}) combined with Theorems 431 and 432 of \cite{hardy}, which state that the normal orders of $\omega(u)$ and $\tau(u)$ are $\log \log u$ and $2^{\log \log u}$, respectively. 

Next, write $p_n$ for the $n$-th prime congruent to 1 mod 4. It is a well known fact that
\begin{equation}
\label{n_prime}
p_n = O(n \log n).
\end{equation}
For each $n \geq 1$, define $v_n = \prod_{i=1}^n p_i^2$. Write $\N_I(v_n)$ for the number of lattices in $\WR(\zed^2)$ with determinant $v_n$ that come from ideals in $\zed[i]$, then
\begin{equation}
\label{NI}
\N(v_n) \geq \N_I(v_n) = 2\alpha(v_n)+1 = 3^n.
\end{equation}

%%%%%%%%%%%%%%%%%%%%%%%%%%%%%%%%%%%%%%%%%%%%%%%%%%%%%%%%%%%%%%%%%%
%the number of divisors of $v_n$ that are products of squares of $i$ primes is $\binom{n}{i}$ for each $0 \leq i \leq n$, and $\alpha_*$ evaluated at each such divisor is $2^{i-1}$. Then, by Theorem \ref{count}
%$$\N(v_n) \geq \beta(v_n) + 2 \mathop{\sum_{b|v_n, 1 \leq b<v_n}}_{b\ \text{a square}} \alpha_* \left(\frac{v_n}{b} \right) (2\beta(b)-1) \geq \sum_{i=0}^n 2^i \binom{n}{i} = 3^n.$$
%%%%%%%%%%%%%%%%%%%%%%%%%%%%%%%%%%%%%%%%%%%%%%%%%%%%%%%%%%%%%%%%%%

\noindent
Let $k$ be a positive integer. We want to choose $n$ such that
\begin{equation}
\label{nO_1}
\N(v_n) \geq 3^n \geq \left( \log \left( \prod_{i=1}^n p_i^2 \right) \right)^k = 2^k \left( \sum_{i=1}^n \log p_i \right)^k.
\end{equation}
By (\ref{n_prime}),
\begin{equation}
\label{nO_2}
\sum_{i=1}^n \log p_i = \sum_{i=1}^n \log \left( O(i \log i) \right) = \sum_{i=1}^n O \left( \log (i \log i) \right) = \sum_{i=1}^n O \left( \log i  \right) \leq O(n \log n).
\end{equation}
Combining (\ref{nO_1}) with (\ref{nO_2}) and taking logarithms, we see that it is sufficient to choose $n$ such that
$$\frac{n}{\log n} \geq O(k),$$
hence we can take $n = O(k \log k)$. Then, by (\ref{nO_2}), for this choice of $n$ we have
$$v_n =  \exp \left( 2 \sum_{i=1}^n \log p_i \right) \leq \exp \left( O(n \log n) \right) =  \exp \left( O (k (\log k)^2) \right).$$
Let $u_k = v_n$ for this choice of $n$, and so $n = \omega(u_k)$. This completes the proof.
\endproof

%%%%%%%%%%%%%%%%%%%%%%%%%%%%%%%%%%%%%%%%%%%%%%%%%%%%%%%%%%%%%%%%%%%%%%%%%%%%
%\begin{center} 
%\begin{tabular}{|l|l|l|l|} \hline
%{\em $u$} & {\em $\N(u)$} & {\em $(\log u)^{\log 8}$} & {\em $\frac{\N(u)}{(\log u)^{\log 8}}$} \\ \hline \hline
%$(43 \times 47)^{40}$ & 18642 & 145994.9689 & 0.12768933 \\ \hline
%$(13 \times 17)^{40}$ & 19946 & 71458.31778 & 0.27912776 \\ \hline
%$(43 \times 47 \times 53)^{10}$ & 49620 & 19565.41976 & 2.5361071 \\ \hline
%$(61 \times 73 \times 89 \times 97)^6$ & 151414 & 15889.87537 & 9.5289608 \\ \hline
%$(43 \times 47 \times 53 \times 61 \times 73 \times 89)^2$ & 39437 & 3262.998184 & 12.086124 \\ \hline
%$(43 \times 47 \times 53 \times 61)^6$ & 203364 & 12720.16412 & 15.98753 \\ \hline
%$(43 \times 47 \times 53 \times 61 \times 73)^4$ & 393868 & 9048.887758 & 43.526675 \\ \hline
%$(43 \times 47 \times 53 \times 61 \times 73 \times 89 \times 97)^2$ & 273738 & 4660.008800 & 58.741949 \\ \hline
%$(13 \times 17 \times 53 \times 61 \times 73 \times 89 \times 97 \times 101)^2$ & 1429871 & 5497.165157 & 260.11061 \\ \hline
%$(43 \times 47 \times 53 \times 61 \times 73 \times 89 \times 97 \times 101)^2$ & 2052339 & 6332.241315 & 324.10941 \\ \hline
%$(43 \times 47 \times 53 \times 61 \times 73 \times 89)^4$ & 6487069 & 13790.85240 & 470.38927 \\ \hline
%$(5 \times 13 \times 17 \times 29 \times 37 \times 43 \times 47 \times 53)^2$ & 1743637 & 3574.193896 & 487.84063 \\ \hline
%\end{tabular}
%\end{center}
%%%%%%%%%%%%%%%%%%%%%%%%%%%%%%%%%%%%%%%%%%%%%%%%%%%%%%%%%%%%%%%%%%%%%%%%%%%%

\begin{rem} Let $v_n$ be as in the proof of Corollary \ref{size_N} above, i.e. $v_n = \prod_{i=1}^n p_i^2$, where $p_1, p_2, \dots$ are primes congruent to 1 mod 4; for instance, the first 9 such primes are 5, 13, 17, 29, 37, 43, 47, 53, 61. For each $k$ choose the smallest $n$ so that $v_n > (\log v_n)^k$, and let $u_k = v_n$ for this choice of $n$. Here is the actual data table for the first few values of the sequence $\{u_k\}$ computed with Maple.

\begin{center} 
\begin{tabular}{|l|l|l|l|l|} \hline
{\em $k$} & {\em $n$} & {\em $u_k = v_n$} & {\em $\N(u_k)$} & {\em $(\log u_k)^k$} \\ \hline \hline
1 & 2 & $4225$ & 9 & 8.34877454 \\ \hline
2 & 4 & $1026882025$ & 518 & 430.5539044 \\ \hline
3 & 7 & $5741913252704971225$ & 215002 & 80589.79464 \\ \hline
4 & 9 & $60016136730202390980384025$ & 14324372 & 12413026.85 \\ \hline
\end{tabular}
\end{center}

\noindent
Notice that the choice of $n = O(k \log k)$ as in Corollary \ref{size_N} insures that not just $\N(u_k)$, but even the much smaller $\N_I(u_k)$ (compare for instance the values of $\N(u_k)$ in the table above to $\N_I(u_k) = 3^n$) is greater than $(\log u_k)^k$, and even with this stronger restriction $u_k$ and $\omega(u_k)$ grow relatively slow as functions of $k$.
\end{rem}
\bigskip

\section{Counting well-rounded lattices with fixed minimum}

Let $m \in \zed_{>0}$, then, as stated in \cite{martinet:venkov}, there exist $\left[ \frac{m+1}{2} \right]$ WR lattices $\Lambda$, not necessarily integral, of rank 2 in $\real^2$ with $|\Lambda| = \sqrt{m}$, generated by a minimal basis $\bx,\bwy$ with $0 < \bx^t \bwy \leq \left[ \frac{m-1}{2} \right]$. This information, however, does not lead to an explicit formula for the number of WR sublattices of $\zed^2$ of prescribed minimum. We derive such a formula here.
\smallskip

Let $m \in \zed_{>0}$. Suppose that $\Lambda \in \WR(\zed^2)$ and $|\Lambda|^2=m$, then by Lemma \ref{gauss} there exists a representation $\Lambda = \left( \begin{matrix} x_1&y_1 \\ x_2&y_2 \end{matrix} \right) \zed^2$ with 
\begin{equation}
\label{abcd_min}
\bx = \left(\begin{matrix} x_1 \\ x_2 \end{matrix} \right),\  \bwy = \left(\begin{matrix} y_1 \\ y_2 \end{matrix} \right) \in \zed^2,\ \|\bx\|=\|\bwy\|=m,\ \theta(\bx,\bwy) \in \left[ \frac{\pi}{3}, \frac{\pi}{2} \right],
\end{equation}
where $\theta(\bx,\bwy)$ is the angle between vectors $\bx$ and $\bwy$, since if $\theta(\bx,\bwy) \in \left( \frac{\pi}{2}, \frac{2\pi}{3} \right]$ we can always replace $\bx$ with $-\bx$ or $\bwy$ with $-\bwy$ to ensure that $\theta(\bx,\bwy) \in \left[ \frac{\pi}{3}, \frac{\pi}{2} \right]$. Then define
$$\C_m = \{ \bx \in \zed^2 : x_2 > 0,\ \|\bx\| = m \}.$$
For each $\bx \in \C_m$ let
$$E_m(\bx) = \left\{ \bwy \in \zed^2 : y_2 \geq 0,\ \|\bwy\| = m,\ \theta(\bx,\bwy) \in \left[ \frac{\pi}{3}, \frac{\pi}{2} \right] \right\},$$
and define $\eta_m(\bx) = |E_m(\bx)|$. The following result follows immediately.

\begin{thm} \label{count_min} Let $m \in \Mm$. Let $\N'(m)$ be the number of lattices in $\WR(\zed^2)$ with minimum equal to $m$. Then
$$\N'(m) = \sum_{\bx \in \C_m} \eta_m(\bx).$$
\end{thm}

Notice that for each $\bx \in \C_m$, $\eta_m(\bx)$ is precisely the number of integer lattice points on the arc of the circle of radius $\sqrt{m}$, bounded by the points $\left( \begin{matrix} \frac{x_1}{2} - \frac{x_2\sqrt{3}}{2} \\ \frac{x_1\sqrt{3}}{2} + \frac{x_2}{2} \end{matrix} \right)$ and $\left( \begin{matrix} -x_2 \\ x_1 \end{matrix} \right)$. The angle corresponding to this arc is $\frac{\pi}{6}$. On the other hand, $\alpha(m)$ as defined by (\ref{alpha}) is the number of integer lattice points on the quarter-circle of radius $\sqrt{m}$ centered at the origin and bounded by the points $(\sqrt{m},0)$, $(0, \sqrt{m})$. It is not difficult to see that for every $\bx \in \C_m$, 
\[ \eta_m(\bx) \leq \left\{ \begin{array}{ll}
\alpha(m) & \mbox{if $m$ is not a square} \\
\alpha(m)+1 & \mbox{if $m$ is a square.}
\end{array}
\right. \]
Indeed, for each $\left( \begin{matrix} 0 \\ \sqrt{m} \end{matrix} \right) \neq \left( \begin{matrix} y_1 \\ y_2 \end{matrix} \right) \in E_m(\bx)$, either $\left( \begin{matrix} y_1 \\ y_2 \end{matrix} \right)$ or $\left( \begin{matrix} -y_2 \\ y_1 \end{matrix} \right)$ is contained in the first quadrant and so is counted by $\alpha(m)$; if $m$ is a square, we add one to account for the point $\left( \begin{matrix} 0 \\ \sqrt{m} \end{matrix} \right)$, which is not counted by $\alpha(m)$. In general, these bounds are sharp, for instance $\alpha(13)=\eta_{13}\left( \begin{matrix} 2 \\ 3 \end{matrix} \right)=1$. However, if $\alpha(m)$ is large and the integral points $\bx$ are well distributed on the quarter-circle, it is possible to do better. For this $m$ needs to satisfy certain special conditions. More precisely, write prime decomposition of $m$ as
$$m = 2^w p_1^{2l_1} \dots p_s^{2l_s} q_1^{k_1} \dots q_r^{k_r},$$
where $p_i \equiv 3\ (\md 4)$, $q_j \equiv 1\ (\md 4)$, $w \in \zed_{\geq 0}$, $l_i \in \frac{1}{2} \zed_{>0}$, and $k_j \in \zed_{>0}$ for all $1 \leq i \leq s$, $1 \leq j \leq r$. If $l_i \notin \zed$ for any $1 \leq i \leq s$, then $\alpha(m) = 0$, so let us assume that $l_i \in \zed$ for all $1 \leq i \leq s$. Then
$$\alpha(m) = \alpha \left( \frac{m}{2^w p_1^{2l_1} \dots p_s^{2l_s}} \right),$$
hence we can assume that
\begin{equation}
\label{m_cong1}
m = q_1^{k_1} \dots q_r^{k_r},
\end{equation}
where $q_j \equiv 1\ (\md 4)$ and $k_j \in \zed_{>0}$ for all $1 \leq j \leq r$. Define
$$L(m) = \sqrt{ \frac{\log (q_1 \dots q_r)}{\log \alpha(m)} },$$
and let 
$$\Mm_1 = \left\{ m \in \Mm : m\ \text{as in (\ref{m_cong1}) and } L(m) \rightarrow 0 \text{ as } m \rightarrow \infty \right\}.$$
A result of Babaev \cite{babaev} implies that for $m \in \Mm_1$
\begin{equation}
\label{alphap_bound}
\eta_m(\bx) = \frac{\alpha(m)}{3} + O(L(m)\alpha(m)),
\end{equation}
for each $\bx \in \C_m$. For $m \in \Mm \setminus \Mm_1$ I am not aware of upper bounds on $\eta_m(\bx)$ better than $\alpha(m)$; a classical result of Jarnik on the number of integral lattice points on convex curves \cite{jarnik} as well as more modern results, for instance of Bombieri and Pila \cite{pila}, imply a general bound on $\eta_m(\bx)$ which is at best $O \left(m^{\frac{1}{4}+\eps} \right)$ for each $\bx \in \C_m$.
\bigskip

In precisely the same manner as Corollary \ref{gen} follows from Theorem \ref{count}, the following is an immediate consequence of Theorem \ref{count_min}.

\begin{cor} \label{gen_min} Let $A \in O_2(\real)$. Then for each $m \in \Mm$ the number of full-rank WR sublattices of $A\zed^2$ with squared minimum equal to $m$ is given by $\N'(m)$ as in Theorem \ref{count_min}.
\end{cor}
\bigskip

\begin{rem} \label{min_det} Notice that even fixing both, the minimum and the determinant, does not identify an element of $\WR(\zed^2)$ uniquely. For instance, if $u$ is representable as a sum of two squares, then the number of lattices $\Lambda \in \WR(\zed^2)$ with $|\Lambda|^2 = \det(\Lambda) = u$ is
\[ \N_I(u) = \left\{ \begin{array}{ll}
2\alpha(u) & \mbox{if $u$ is not a square} \\
2\alpha(u)+1 & \mbox{if $u$ is a square,}
\end{array}
\right. \]
i.e. precisely the number of lattices in $\WR(\zed^2)$ of determinant $u$ coming from ideals in $\zed[i]$, as defined in section~6. Hence even this number can tend to infinity with~$u$. 
\end{rem}
\bigskip

\section{Zeta function of well-rounded lattices}

Given any finitely generated group $G$, it is possible to associate a zeta function $\zeta_G(s) = \sum_{n=1}^{\infty} a_n n^{-s}$ to it, where the coefficients $a_n$ count the number of its subgroups of index $n$ and $s \in \cee$ (see \cite{lubot}, Chapter 15 for details). Such zeta functions are extensively studied objects, since they encode important arithmetic information about the group in question and often have interesting properties. For example, by Theorem 15.1 of \cite{lubot} (see also (5) of \cite{reiner})
\begin{equation}
\label{all_zeta}
\zeta_{\zed^2}(s) = \sum_{\Lambda \subseteq \zed^2} (\det(\Lambda))^{-s} = \sum_{u=1}^{\infty} F(2,u) u^{-s} = \zeta(s) \zeta(s-1),
\end{equation}
where the sum is taken over all sublattices $\Lambda$ of $\zed^2$ of finite index, $F(2,u)$ is given by (\ref{all_sublattices}), and $\zeta(s)$ is the Riemann zeta function. This $\zeta_{\zed^2}(s)$ is an example of Solomon's zeta function (see \cite{reiner}, \cite{solomon}). The identity (\ref{all_zeta}) holds in the half-plane $\Re(s) > 2$, where this series is absolutely convergent, and so the function $\zeta_{\zed^2}(s)$ is analytic (Proposition 1 of \cite{reiner}). Moreover, in this half-plane $\zeta_{\zed^2}(s) = \sum_{u=1}^{\infty} \sigma(u) u^{-s}$, where $\sigma(u) = \sum_{d|u} d$ (see Theorem 290 of \cite{hardy}); $\zeta_{\zed^2}(s)$ has a pole at $s=2$.

In this section we study the properties of the partial zeta function corresponding not to all, but only to the {\it well-rounded} sublattices of $\zed^2$ as defined by (\ref{WR_zeta}). We also define the Dedekind zeta function of Gaussian integers $\zed[i]$
\begin{equation}
\label{dedekind}
\zeta_{\zed[i]}(s) = \sum_{\Aa \subseteq \zed[i]} \Nn(\Aa)^{-s}  = \mathop{\sum_{\Lambda \in \WR(\zed^2)}}_{|\Lambda|^2=\det(\Lambda)} (\det(\Lambda))^{-s} = \sum_{m=1}^{\infty} \N_I(m) m^{-s},
\end{equation}
where $\N_I(m)$ is as above, and the first sum is taken over all the ideals $\Aa = (a+bi)\zed[i]$ for some $a,b \in \zed$, and $\Nn(\Aa) = a^2+b^2$ is the norm of such ideal. In other words, coefficients of $\zeta_{\zed[i]}(s)$ count the elements of $\WR(\zed^2)$ that come from ideals in $\zed[i]$ while coefficients of $\zeta_{\WR(\zed^2)}(s)$ count all elements of $\WR(\zed^2)$. We also note that $\zeta_{\zed[i]}(s)$ is analytic on $\Re(s) > 1/2$ except for a simple pole at $s=1$ (see Theorem 5 on p. 161 of \cite{lang}). It is clear that for all $u \in \zed_{>0}$
\begin{equation}
\label{zeta_c}
\N_I(u) \leq \N(u) \leq F(2,u),
\end{equation}
in other words $\zeta_{\WR(\zed^2)}(s)$ is ``squeezed'' between $\zeta_{\zed[i]}(s)$ and $\zeta_{\zed^2}(s)$. Moreover, our estimates on coefficients in the previous sections suggest that $\zeta_{\WR(\zed^2)}(s)$ should be ``closer'' to $\zeta_{\zed[i]}(s)$ than to $\zeta_{\zed^2}(s)$. Theorem \ref{zeta}, which we will now prove, makes this statement more precise. We start by studying some related Dirichlet series.

\begin{lem} \label{Dir1} Let $t = t(u)$ be as in Theorem \ref{count}. The Dirichlet series $\sum_{u=1}^{\infty} \frac{2\alpha_*(t)}{u^{s}}$ is absolutely convergent at least in the half-plane $\Re(s) > 1$ with a simple pole at $s=1$. Moreover, when $\Re(s) > 1$ it has an Euler product expansion:
\begin{equation}
\label{Euler_prod}
\sum_{u=1}^{\infty} \frac{2\alpha_*(t)}{u^{s}} = \left( 1 + \frac{1}{2^s} + \frac{1}{4^s} \right) \prod_{p \equiv 1 (\md 4)} \frac{p^s+1}{p^s-1}.
\end{equation}
\end{lem}

\proof
First of all, notice that for every $u \in \zed_{>0}$, 
$$2\alpha_*(t) \leq 2\alpha(t) \leq 2\alpha(u) + 1 \leq \N_I(u) + 1,$$
therefore  $\sum_{u=1}^{\infty} 2\alpha_*(t) u^{-s}$ is absolutely convergent at least on the half-plane $\Re(s) > 1$ with at most a simple pole at $s=1$, since $\sum_{u=1}^{\infty} (\N_I(u)+1) u^{-s} = \zeta_{\zed[i]}(s) + \zeta(s)$ is. Next, let
\[ \alpha_*'(n) = \left\{ \begin{array}{ll}
\alpha_*(n) & \mbox{if $n$ is odd} \\
0 & \mbox{if $n$ is even.}
\end{array}
\right. \]
Notice that $2\alpha'_*$ is a multiplicative arithmetic function, specifically $2\alpha'_*(1)=1$ and $2\alpha'_*(mn) = 2\alpha'_*(m) 2\alpha'_*(n)$ for all $m,n \in \zed_{>0}$ with $\gcd(m,n)=1$. Therefore, by Theorem 286 of \cite{hardy} the series $\sum_{u=1}^{\infty} 2\alpha'_*(u) u^{-s}$ has the following Euler-type product representation, where $p$ is always a prime:
\begin{eqnarray*}
\sum_{u=1}^{\infty} 2\alpha'_*(u) u^{-s} & = & \prod_p \left( \sum_{k=0}^{\infty} 2\alpha'_*(p^k) p^{-ks} \right) = \prod_{p \equiv 1 (\md 4)} \left( 1 + 2 \sum_{k=1}^{\infty} p^{-ks} \right) \\
& = & \prod_{p \equiv 1 (\md 4)} \left( \frac{2}{1-p^{-s}} - 1 \right) = \prod_{p \equiv 1 (\md 4)} \frac{p^s+1}{p^s-1},
\end{eqnarray*}
whenever this product is convergent. Also notice that since $\alpha_*(2u) = 0$ if $2|u$ and $\alpha_*(2u) = \alpha_*(u)$ if $2 \nmid u$, we have
\begin{eqnarray*}
\sum_{u=1}^{\infty} \frac{2\alpha_*(t)}{u^{s}} & = & \sum_{u=1}^{\infty} \frac{2\alpha'_*(u)}{u^{s}} + \sum_{u=1}^{\infty} \frac{2\alpha_*(u)}{(2u)^{s}} \\
 & = & \sum_{u=1}^{\infty} \frac{2\alpha'_*(u)}{u^{s}} + \frac{1}{2^s} \left( \sum_{u=1}^{\infty} \frac{2\alpha'_*(u)}{u^{s}} +  \frac{1}{2^s} \sum_{u=1}^{\infty} \frac{2\alpha'_*(u)}{u^{s}} \right) \\
& = & \left( 1 + \frac{1}{2^s} + \frac{1}{4^s} \right) \sum_{u=1}^{\infty} \frac{2\alpha'_*(u)}{u^{s}},
\end{eqnarray*}
which proves (\ref{Euler_prod}) when $\prod_{p \equiv 1 (\md 4)} \frac{p^s+1}{p^s-1}$ is convergent. It is easy to notice that this happens when $\Re(s) > 1$, but $\prod_{p \equiv 1 (\md 4)} \frac{p+1}{p-1}$ diverges, meaning that $\sum_{u=1}^{\infty} \frac{2\alpha_*(t)}{u^{s}}$ must have a pole at $s=1$, and by our argument above we know that it must be a simple pole. This completes the proof.
\endproof

For the next lemma, let $\B_{\nu}$ be as in (\ref{prod_set}) in section 4.

\begin{lem} \label{Dir2} For each $1 < \nu \leq 3^{1/4}$, the Dirichlet series $\sum_{u \in \B_{\nu}} \frac{1}{u^{s}}$ is absolutely convergent in the half-plane $\Re(s) > 1$ with a simple pole at $s=1$ in the sense of~(\ref{pole_def}).
\end{lem}

\proof
First notice that 
$$\sum_{u \in \B_{\nu}} \left| \frac{1}{u^s} \right| \leq \sum_{n=1}^{\infty} \left| \frac{1}{n^s} \right|,$$
and so must be analytic when $\Re(s) > 1$ with at most a simple pole at $s=1$. 
\smallskip

On the other hand, let the Dirichlet lower density of the set $\B_{\nu}$ be defined as
$$\liminf_{s \rightarrow 1^+} \frac{\sum_{u \in \B_{\nu}} u^{-s}}{\sum_{u \in \zed_{>0}} u^{-s}} = \liminf_{s \rightarrow 1^+} \frac{1}{\zeta(s)} \sum_{u \in \B_{\nu}} u^{-s}.$$
It is a well known fact (see for instance equation (1.6) of \cite{ahlswede}) that the Dirichlet lower density of a set is greater or equal than its lower density. Hence, by Theorem~\ref{B_density}
$$0 < \frac{\nu-1}{2 \nu} \leq \DL_{\B_{\nu}} \leq \liminf_{s \rightarrow 1^+} \frac{1}{\zeta(s)} \sum_{u \in \B_{\nu}} u^{-s},$$
which implies that $\sum_{u \in \B_{\nu}} u^{-s}$ must have a pole of the same order as $\zeta(s)$ at $s=1$. This completes the proof.
\endproof

\begin{lem} \label{Dir3}  For each $1 < \nu \leq 3^{1/4}$, Dirichlet series $\sum_{u=1}^{\infty} \frac{\beta_{\nu}(u)}{u^{s}}$ is absolutely convergent in the half-plane $\Re(s) > 1$, and is bounded below by a Dirichlet series with a pole of order 1 at $s=1$. Moreover, for every real $\eps >0$ there exists a Dirichlet series with a pole of order $1+\eps$ at $s=1$, which bounds $\sum_{u=1}^{\infty} \frac{\beta_{\nu}(u)}{u^{s}}$ from above.
\end{lem}

\proof
For each $1 < \nu \leq 3^{1/4}$, define $\chi_{\nu}$ to be the characteristic function of the set $\B_{\nu}$, i.e. for each $u \in \zed_{>0}$,
\[ \chi_{\nu}(u) = \left\{ \begin{array}{ll}
1 & \mbox{if $u \in \B_{\nu}$} \\
0 & \mbox{if $u \notin \B_{\nu}$.}
\end{array}
\right. \]
Clearly, $\beta_{\nu}(u) \geq \chi_{\nu}(u)$, therefore
$$\sum_{u=1}^{\infty} \left| \frac{\beta_{\nu}(u)}{u^{s}} \right| \geq \sum_{u=1}^{\infty} \left| \frac{\chi_{\nu}(u)}{u^{s}} \right| = \sum_{u \in \B_{\nu}} \left| \frac{1}{u^{s}} \right|,$$
which, combined with Lemma \ref{Dir2}, proves the lower bound of the lemma.
\smallskip

On the other hand, recall that $\beta_{\nu}(u) \leq \Delta(u)$ for all $u \in \zed_{>0}$, where $\Delta(u)$ is Hooley's $\Delta$-function, as defined in section~2, hence $\sum_{u=1}^{\infty} \left| \beta_{\nu}(u) u^{-s} \right| \leq \sum_{u=1}^{\infty} \left| \Delta(u) u^{-s} \right|$. Hooley's $\Delta$-function is known to satisfy
\begin{equation}
\label{hooley_log}
\sum_{u=1}^{\infty} \left| \frac{\Delta(u)}{u^{s}} \right| \ll_{\eps} \sum_{u=1}^{\infty} \left| \frac{(\log u)^{\eps}}{u^{s}} \right|,
\end{equation}
for every $\eps > 0$, which is a consequence of Tenenbaum's bound on the average order of $\Delta(u)$ (see \cite{tenenbaum}, also \cite{hall}), and so the upper bound of the lemma follows by observing that $\sum_{u=1}^{\infty} (\log u)^{\eps} u^{-s}$ has a pole of order $1+\eps$ at $s=1$. Since $\sum_{u=1}^{\infty} (\log u)^{\eps} u^{-s}$ is absolutely convergent in the half-plane $\Re(s) > 1$, (\ref{hooley_log}) also proves that so is $\sum_{u=1}^{\infty} \frac{\beta_{\nu}(u)}{u^{s}}$.
\endproof

We are now ready to prove Theorem \ref{zeta}.
\bigskip

\noindent
{\it Proof of Theorem \ref{zeta}.}
First of all notice that (\ref{zeta_c}) combined with comparison test for series imply that $\sum_{u=1}^{\infty} \N(u) u^{-s}$ has a pole at $s=1$, since $\zeta_{\zed[i]}(s)$ has a pole at $s=1$, and is absolutely convergent, i.e. $\zeta_{\WR(\zed^2)}(s)$ is analytic, for $\Re(s) > 2$, since $\zeta_{\zed^2}(s)$ is analytic when $\Re(s) > 2$. In fact, we can do better. Let
\[ \beta'(n) = \left\{ \begin{array}{ll}
2\beta(n)-1 & \mbox{if $n$ is a square} \\
2\beta(n) & \mbox{if $n$ is not a square,}
\end{array}
\right. \]
for every $n \in \zed_{>0}$. Notice that for every $u \in \zed_{>0}$, $\N(u)$ can be expressed in terms of the Dirichlet convolution of arithmetic functions $2\alpha_*$ and $\beta'$:
$$\N(u) = (2\alpha_* * \beta')(t) + \left( \delta_1(t)\beta(t) + \delta_2(t)\beta\left( \frac{t}{2} \right) - 2\beta'(t) - \frac{1+(-1)^t}{2} \beta' \left( \frac{t}{2} \right) \right),$$
where $t=t(u)$, $\delta_1(t)$, and $\delta_2(t)$ are as in Theorem \ref{count}. Therefore, by Theorem 284 of \cite{hardy}
\begin{eqnarray}
\label{zeta_prod1}
\zeta_{\WR(\zed^2)}(s) & = & \sum_{u=1}^{\infty} (2\alpha_* * \beta')(t) u^{-s} \nonumber \\
& + & \sum_{u=1}^{\infty}  \left( \delta_1(t)\beta(t) + \delta_2(t)\beta\left( \frac{t}{2} \right) - 2\beta'(t) - \frac{1+(-1)^t}{2} \beta' \left( \frac{t}{2} \right) \right) u^{-s} \nonumber \\
& = & \left( \sum_{u=1}^{\infty} 2\alpha_*(t) u^{-s} - 2\right) \left( \sum_{u=1}^{\infty} \beta'(t) u^{-s} \right) \nonumber \\
& + & \sum_{u=1}^{\infty}  \left( \delta_1(t)\beta(t) + \delta_2(t)\beta\left( \frac{t}{2} \right) \right) u^{-s} - \sum_{u=1}^{\infty} \frac{1+(-1)^t}{2} \beta' \left( \frac{t}{2} \right) u^{-s},
\end{eqnarray}
whenever these three series are absolutely convergent. Now notice that 
$$\frac{1}{|2^s|} \sum_{u=1}^{\infty} \frac{\beta(u)}{|u^s|} \leq \sum_{u=1}^{\infty} \frac{\delta_1(t)\beta(t)}{|u^s|} \leq 2 \sum_{u=1}^{\infty} \frac{\beta(u)}{|u^s|},$$
and
$$\frac{1}{|4^s|} \sum_{u=1}^{\infty} \frac{\beta(u)}{|u^s|} \leq \sum_{u=1}^{\infty} \frac{\delta_2(t)\beta\left(\frac{t}{2}\right)}{|u^s|} = \frac{1}{|4^s|} \sum_{u=1}^{\infty} \frac{\delta_1(u)\beta(u)}{|u^s|} \leq \frac{2}{|4^s|} \sum_{u=1}^{\infty} \frac{\beta(u)}{|u^s|},$$ 
as well as
$$\frac{1}{|2^s|} O \left( \sum_{u=1}^{\infty} \frac{\beta(u)}{|u^s|} \right) = \frac{1}{|2^s|} \sum_{u=1}^{\infty} \frac{\beta'(u)}{|u^s|} \leq \sum_{u=1}^{\infty} \frac{\beta'(t)}{|u^s|} \leq \sum_{u=1}^{\infty} \frac{\beta'(u)}{|u^s|} = O \left( \sum_{u=1}^{\infty} \frac{\beta(u)}{|u^s|} \right),$$
whenever $\Re(s) > 0$. Also
$$\sum_{u=1}^{\infty} \frac{1+(-1)^t}{2} \beta' \left( \frac{t}{2} \right) u^{-s} = 4^{-s} \sum_{u=1}^{\infty} \beta'(u) u^{-s}.$$
Now the conclusion of Theorem \ref{zeta} follows by applying these observations along with Lemmas \ref{Dir1}, \ref{Dir3} to (\ref{zeta_prod1}).
\boxed{ }
\bigskip

\begin{rem} \label{ded_zeta_sq} Notice that one implication of Theorem \ref{zeta} is that $\N(u)$, the coefficient of $\zeta_{\WR(\zed^2)}(s)$, grows, roughly speaking, like the coefficient of $\zeta_{\zed[i]}(s)^2$, which~is
$$\sum_{mn = u} \N_I(m) \N_I(n).$$
\end{rem}

Finally, we mention that in the same manner one can define zeta functions of well-rounded sublattices of any lattice $\Omega$ in $\real^N$ for any $N$. Studying the properties of these functions may yield interesting arithmetic information about the distribution of such sublattices.
\bigskip

{\bf Acknowledgment.} I would like to thank Kevin Ford, Preda Mihailescu, Baruch Moroz, Gabriele Nebe, Bogdan Petrenko, Sinai Robins, Eugenia Soboleva, Paula Tretkoff, Jeff Vaaler and Victor Vuletescu for their helpful comments on the subject of this paper. I would also like to acknowledge the wonderful hospitality of Max-Planck-Institut f\"{u}r Mathematik in Bonn, Germany, where a large part of this work has been done.

%\nocite{*}
\bibliographystyle{plain}  % Here the bibliography 
\bibliography{esm}   % is inserted.

\end{document}